\documentclass{article}
\usepackage[english]{babel}

\usepackage[letterpaper,top=2cm,bottom=2cm,left=2cm,right=2cm,marginparwidth=1.75cm]{geometry}

\usepackage{calrsfs}
\usepackage{amssymb}
\usepackage{amsmath}
\usepackage{amsthm}
\usepackage{subcaption}
\usepackage{multirow}
\usepackage{graphicx}
\usepackage{amsfonts}
\usepackage{mathtools}
\usepackage{mathrsfs}
\usepackage{placeins}
\usepackage{textgreek}
\usepackage{booktabs}
\newtheorem{remark}{Remark}

\usepackage[table]{xcolor}
\usepackage{booktabs}
\usepackage{longtable}
\usepackage{array}

\usepackage[mathscr]{euscript}
\usepackage[dvipsnames]{xcolor}
\usepackage[colorlinks=true, allcolors=OliveGreen]{hyperref}
\usepackage{filecontents}
\usepackage{tikz}
\usepackage{pgfplots}
\usepackage{authblk}
\usepackage{todonotes}
\usepackage{array}
\newcolumntype{C}{>{\centering\arraybackslash}p{2.5cm}}

\usepackage[dvipsnames,x11names]{xcolor}

\usepackage{graphicx} 
\pgfplotsset{compat=1.18}
\usepackage{caption}
\usepackage{subcaption}
\usepackage{float}
\usepackage{graphics} 
\usepackage{tikz}
\usepackage{csvsimple}
\usepackage{siunitx}
\usepackage{comment}
\usepackage{tikz} 
\usepackage{pgfmath,pgffor,pgfplots}
\pgfplotsset{compat=1.18} 
\usetikzlibrary{external}
\usetikzlibrary{pgfplots.groupplots}
\usepackage{orcidlink}
\usepackage[utf8]{inputenc}
\usetikzlibrary{arrows.meta,decorations.pathmorphing,backgrounds,positioning,fit,petri}

\newcommand{\averagel}{\{\!\!\{}
\newcommand{\averager}{\}\!\!\}}

\newcommand{\jumpl}{[\![}
\newcommand{\jumpr}{]\!]}

\DeclareMathAlphabet{\mathcalligra}{T1}{calligra}{m}{n}

\newcommand{\Napl}{\mathrm{Na}^+}
\newcommand{\Capl}{\mathrm{Ca}^{2+}}
\newcommand{\Kpl}{\mathrm{K}^+}
\newcommand{\Na}{[\mathrm{Na}^+]_i}
\newcommand{\Ca}{[\mathrm{Ca}^{2+}]_i}
\newcommand{\K}{[\mathrm{K}^+]_o}
\newcommand{\Kbath}{\mathrm{K}_\mathrm{bath}}
\newcommand{\Abeta}[1]{A\textbeta{#1}}
\newcommand{\OmegaAb}{\Omega_{\mathrm{A}\beta}}
\newcommand{\Ab}{[A\mathrm{\beta}]}
\newcommand{\muM}{\mu\mathrm{M}}

\graphicspath{{Photos/}}

\tikzset{  font={\fontsize{15pt}{12}\selectfont}}

\title{A novel mathematical and computational framework of amyloid-beta triggered seizure dynamics in Alzheimer's disease\footnote{\textbf{Funding}: This work is partially funded by the European Union (ERC SyG, NEMESIS, project number 101115663). Views and opinions expressed are, however, those of the authors only and do not necessarily reflect those of the European Union or the European Research Council Executive Agency. Neither the European Union nor the granting authority can be held responsible for them. CBLS has been funded by the National Recovery and Resilience Plan (NRRP), Mission 4, Component 1 – Investment 3.4 and Investment 4.1, funded by the European Union. PFA has been partially supported by ICSC—Centro Nazionale di Ricerca in High Performance Computing, Big Data, and Quantum Computing, funded by the European Union-NextGeneration EU.  The present research is part of the activities of the Dipartimento di Eccellenza 2023-2027 grant, funded by MUR. 
The authors acknowledge the CINECA award under the ISCRA initiative, for the availability of high-performance computing resources and support under the project IsCc9\_NeuroDG, PI M. Corti, 2025–2026. PFA, SP, MC, and CBLS are members of INdAM-GNCS. }}

\author[1]{Caterina B. Leimer Saglio \orcidlink{0009-0007-7887-919X}}

\affil[1]{MOX-Dipartimento di Matematica, Politecnico di Milano, Piazza Leonardo da Vinci 32, Milan, 20133, Italy}

\author[1]{Mattia Corti \orcidlink{0000-0002-7014-972X}}

\author[1]{Stefano Pagani \orcidlink{0000-0002-6662-3433}}

\author[1]{Paola F. Antonietti \orcidlink{0000-0002-2138-3878}}

\begin{document}
\maketitle

\begin{abstract}
The association of epileptic activity and Alzheimer's disease (AD) has been increasingly reported in both clinical and experimental studies, suggesting that amyloid-\textbeta{} accumulation may directly affect neuronal excitability. Capturing these interactions requires a quantitative description that bridges the molecular alterations of AD with the fast electrophysiological dynamics of epilepsy. We introduce a novel mathematical model that extends the Barreto-Cressman ionic formulation by incorporating multiple mechanisms of calcium dysregulation induced by amyloid-\textbeta{}, including formation of $\Capl$-permeable pores, overactivation of voltage-gated $\Capl$ channels, and suppression of $\Capl$-sensitive potassium currents. The resulting ionic model is coupled with the monodomain equation and discretized using a $p$-adaptive discontinuous Galerkin method on polytopal meshes, providing an effective balance between efficiency and accuracy in capturing the sharp spatiotemporal electrical wavefronts associated with epileptiform discharges. Numerical simulations performed on idealized and realistic brain geometries demonstrate that progressive amyloid-\textbeta{} accumulation leads to severe alterations in calcium homeostasis, increased neuronal hyperexcitability, and pathological seizure propagation. Specifically, high amyloid-\textbeta{} concentrations produce secondary epileptogenic sources and spatially heterogeneous wavefronts, indicating that biochemical inhomogeneities play a critical role in shaping seizure dynamics. These results illustrate how multiscale modeling provides new mechanistic insights into the interplay between neurodegeneration and epilepsy in Alzheimer's disease.
\end{abstract}

\section{Introduction}
\label{sec:introduction}
Alzheimer's disease and epilepsy are neurological disorders that are often associated in older adults, with incidence rates increasing with age. AD is a neurodegenerative disease characterized by progressive memory loss and cognitive decline. It is caused by the abnormal accumulation, aggregation, and spreading of amyloid-\textbeta{} (\Abeta{}) in the central nervous system. It is generated from the amyloid precursor protein (APP) proteolytic cleavage due to the action of specific enzymes \cite{bloom_amyloid_2014}. Additionally, \Abeta{} peptides are known to impair synaptic transmission and disrupt ionic homeostasis \cite{kuchibhotla_abeta_2008}.  Patients affected by this pathology are 5 to 10 times more likely to develop epileptic seizures than healthy individuals of the same age \cite{zhang_clinical_2022}. Epilepsy is a neurological disorder characterized by a persistent predisposition to generate epileptic seizures, resulting from abnormal electrical activity in the brain \cite{beghi_epidemiology_2020}. During the early stages of AD, the risk of an epileptic seizure is sensibly high; moreover, the risk in people with an early-onset is twice that of late-onset AD \cite{kuchibhotla_abeta_2008}. 
\par
From a microscopic perspective, high concentrations of \Abeta{} trigger neurodegenerative processes by forming toxic aggregates, known as plaques. In particular, \Abeta{} oligomers interfere with calcium ($\Capl$) dynamics by inhibiting plasma membrane calcium ATPases (PMCAs), promoting the formation of calcium-permeable membrane pores, overactivating voltage-gated calcium channels (VGCCs), and suppressing calcium-sensitive potassium currents \cite{good_effect_1996,ishii_amyloid_2019}. Because of the relationship between $\Capl$-dynamics and membrane electrical activity, the perturbations usually affect neuronal excitability \cite{steinlein_calcium_2014}, contributing to an altered neuronal excitability, often manifesting as both hyperexcitability and hypoexcitability, and increasing seizure susceptibility \cite{romoli_amyloid_2021,zhang_clinical_2022}. 
\par
The investigation into these pathological phenomena benefits from the recent developments in mathematical modeling of electrophysiology. Indeed, the coordinated activity of voltage-gated ion channels and active transport mechanisms can be accurately described through biophysical ionic models \cite{hodgkin_components_1952,hodgkin_components_1952,barreto_ion_2011,cressman_influence_2009}. Mathematically, the ionic model is a system of Ordinary Differential Equations (ODEs) that connects the voltage-gated channels and the ionic concentrations with the membrane voltage (action potential) of the cells. Unless these models have been extensively used to investigate both physiological \cite{hodgkin_quantitative_1952} and pathological \cite{barreto_ion_2011,cressman_influence_2009} conditions, their application to neurodegenerative conditions, such as Alzheimer's disease, is a developing area of research \cite{latulippe_mathematical_2018,vonbonhorst_impact_2022}.
\par
In this work, we propose a novel computational model that integrates the pathological influence of \Abeta{} accumulation on neuronal electric activity in epileptic conditions. As a starting point, we consider the Barreto-Cressman (BC) ionic model \cite{barreto_ion_2011,cressman_influence_2009}, which describes the electrical activity in neurons in physiological and pathological conditions. 
Our model introduces some modifications to account for the concentration of \Abeta{} mediated by calcium dysregulation \cite{latulippe_mathematical_2018,vonbonhorst_impact_2022}, described in biological studies. Specifically, we incorporate in the Barreto-Cressman ionic model the effects of \Abeta{} on calcium dynamics, including disruption of calcium homeostasis due to the effects of \Abeta{}-mediated changes on plasma membrane calcium ATPases (PMCA) pump \cite{berrocal_altered_2009,berrocal_calmodulin_2012}, calcium-permeable pores \cite{mirdha_aggregation_2024}, L-type voltage-gated calcium channels \cite{kim_effects_2011}, fast-inactivating potassium channels \cite{good_effect_1996}, and $\Capl$-sensitive potassium channels \cite{yamamoto_suppression_2011,yamamoto_amyloid_2021,jhanandas_cellular_2001}. 
The resulting model captures a range of effects induced by \Abeta{} peptides on ion channel kinetics, calcium buffering, and membrane currents. By analyzing different \Abeta{} concentrations corresponding to progressive stages of AD, our framework enables simulation of AP propagation under neurodegenerative conditions. Our model provides a valuable tool for studying the interplay between \Abeta{} peptides and neuronal excitability, as a first attempt to describe the mechanistic links between calcium imbalance, excitotoxicity, and seizures in AD.
\par
To extend the analysis from the single-neuron evolution to the tissue-level scale, we couple a novel ionic model with a monodomain model for the evolution of the transmembrane potential in a space-time domain. The monodomain model is widely used in the context of electrophysiology \cite{potse_comparision_2006,quarteroni_integrated_2017}, and it has been recently proposed to simulate the electrical activity of the brain at the organ-level \cite{schreiner_simulating_2022, leimer_saglio_high-order_2024}. Alternative approaches rely on graph-based network models, where neurons are represented as nodes connected by structural or functional connectivity patterns. While these models are well suited to capture large-scale communication and synchronization phenomena, their computational cost increases dramatically as the number of neurons and synapses grows, making high-resolution, organ-level simulations practically intractable. The monodomain equation, coupled with our ionic model, allows us to simulate seizure propagation in brain tissue under pathological conditions of \Abeta{} presence. This offers insights into how \Abeta{}-induced calcium dysregulation accelerates network-wide hyperexcitability in brain regions near the lesions. 
\par
To assess the quality of the proposed ionic model and validate the impact of the modelling choices, we perform a set of numerical simulations. First, we focus on the 0D ionic model to analyze the effects on the electrical activity of the single neuron. After that, we perform some numerical simulations of the two-dimensional monodomain equation coupled with the proposed ionic model. To discretize the system, we choose a $p-$adaptive Discontinuous Galerkin formulation on polygonal/polyhedral grids (PolyDG) in space \cite{cangiani_hp-version_2017,antonietti_hp_2013} and a semi-implicit second-order discretization in time. The PolyDG method can easily handle the geometric complexity of the brain, due to the possibility of performing mesh agglomeration \cite{antonietti_agglomeration_2022,antonietti_polytopal_2024}. Moreover, it supports high-order approximations, making it highly efficient in the simulation of wave-propagation problems \cite{corti_discontinuous_2023,antonietti_discontinuous_2024}, such as high-frequency electrical activity \cite{leimer_saglio_high-order_2024}. Finally, the choice to use a $p$-adaptive framework \cite{leimer_saglio_p-adaptive_2025} allows for a sensible reduction in computational cost.
\par
We perform some numerical test cases in two-dimensional settings. In particular, we test our model on a realistic brain coronal section, reconstructed from \Abeta{} concentration derived from positron emission tomography (PET). The goal of the simulations is to show that the spatial heterogeneity of \Abeta{} generates a complex interplay between pathological and healthy regions. Areas characterized by high \Abeta{} concentration exhibit self-induced spiking activity, consistent with the autonomous oscillations observed in the parametric sensitivity analyses. The resulting dynamics display asynchronous propagation fronts and multiple activation centers, revealing how biochemical inhomogeneities and axonal anisotropies jointly contribute to the emergence of spatially fragmented seizure patterns in the late stages of the disease.
\par
In Section~\ref{sec:introduction}, we introduce the biological and clinical motivations for this study, highlighting the growing evidence connecting \Abeta{} pathology with increased seizure susceptibility in Alzheimer's patients.  In Section~\ref{sec:model}, we formulate the mathematical model. Starting from the Barreto--Cressman ionic model for epileptic dynamics, we incorporate  \Abeta{-}dependent modifications acting on multiple cellular targets. The resulting system captures the main electrophysiological alterations observed experimentally in  \Abeta{-}affected neurons. In Section~\ref{sec:3}, we introduce the PolyDG discretization for the coupled monodomain-ionic problem, describing the semi-discrete and fully-discrete formulations. The use of a $p$-adaptive method allows accurate and efficient simulations on complex brain geometries, while ensuring numerical stability and scalability for high-order approximations.  
In Section~\ref{sec:0d}, we present a detailed sensitivity analysis of the ionic model with respect to \Abeta{} concentration. In Section~\ref{sec:2D}, we extend the analysis to the tissue level by coupling the modified ionic model with the monodomain equation. Simulations on idealized and realistic two-dimensional brain domains showed that regions with high \Abeta{} concentration become independent epileptogenic sources.
\color{black}
\section{The mathematical model}
\label{sec:model}
In this section, we present the mathematical model to describe the evolution of the action potential inside the brain tissue, taking into account the concentration of agglomerated \Abeta{ }proteins. The derived model couples the monodomain equation with a modified version of the  Barreto–Cressman ionic model that incorporates the \Abeta{-}effects on the ionic concentrations.
\subsection{Neuronal electrophysiology: the monodomain model}
The common choices to efficiently describe spatial transmembrane potential dynamics in the neural tissue are the bidomain and the monodomain models coupled with properly ionic models \cite{schwartz2016analytic,schreiner_simulating_2022}. In this work, we focus on the monodomain model for its efficiency and ability to provide mechanistic insights.  
Given an open, bounded domain $\Omega \in \mathbb{R}^d$, $(d=2,3)$, and a final time $T>0$, we introduce the transmembrane potential $u = u(\boldsymbol{x},t)$ with $u: \Omega \times [0,T] \rightarrow \mathbb{R}$, and the vector $\boldsymbol{y} = \boldsymbol{y}(\boldsymbol{x},t)$ with $\boldsymbol{y}: \Omega \times [0,T] \rightarrow \mathbb{R}^n, n\ge1,$ that contains the ion concentrations and gating variables of the ionic neuronal model. The coupled problem reads: 
\par
for any time $ t \in (0,T]$, find $u=u(\boldsymbol{x},t)$ and $\boldsymbol{y}=\boldsymbol{y}(\boldsymbol{x},t)$ such that:
\begin{subequations}
\label{eq:monodomain}
    \begin{alignat}{3}
    \label{eq:monodomain:pde}
    \chi_m C_m  \frac{\partial u}{\partial t} - \nabla \cdot (\mathbf{\Sigma} \nabla u) +  \chi_m f(u,\boldsymbol{y}) & = 0, & & \quad \mathrm{in} \; \Omega \times (0,T],
    \\[0pt]
    \label{eq:monodomain:ionic}
    \frac{\partial \boldsymbol{y}}{\partial t} + \boldsymbol{m}(u,\boldsymbol{y}) & = \boldsymbol{0}, & & \quad \mathrm{in} \; \Omega \times (0,T],    
    \\[3pt]
    \label{eq:monodomain:bcs}
    \mathbf{\Sigma} \nabla u \cdot \boldsymbol{n} & = 0,  & & \quad\mathrm{on}\;  \partial \Omega  \times (0,T],    
    \\[5pt]
    \label{eq:monodomain:ics}
    u(0) = u^0, \quad \boldsymbol{y}(0) & = \boldsymbol{y}^0,  & & \quad \mathrm{in}\; \Omega.
    \end{alignat}
\end{subequations}
In Equation~\eqref{eq:monodomain:pde}, $\boldsymbol{\Sigma}$ is defined as the conductivity tensor, $\chi_m$ as the membrane capacitance per unit area, $C_m$ as the membrane capacitance, and $f=f(u,\boldsymbol{y})$ represents the action of the ionic currents. Moreover, in Equation~\eqref{eq:monodomain:ionic}, $\boldsymbol{m}=\boldsymbol{m}(u,\boldsymbol{y})$ represents the evolution of the ion concentrations. We impose homogeneous Neumann boundary conditions in Equation~\eqref{eq:monodomain:bcs} with $\boldsymbol{n}$ defined as the normal to the boundary $\partial \Omega$. Finally, we assign the initial conditions $u^0$ and $\boldsymbol{y}^0$ in Equation~\eqref{eq:monodomain:ics}.

\subsection{The Barreto-Cressman ionic model}
\label{subsec:BC}
The electrophysiological behavior of a single neuron in physiological and high-frequency pathological conditions can be described by the Barreto-Cressman (BC) ionic model \cite{cressman_influence_2009}. This is a conductance-based ionic model able to describe the evolution of the neuron’s AP, based on ion dynamics taking into account the main ions responsible for ionic imbalances such as calcium ($\Capl$), potassium ($\Kpl$), and sodium ($\Napl$) concentrations. The Barreto-Cressman ionic model is exploited to analyze and display different neuronal dynamics, i.e., a fast-spiking behavior of the transmembrane potential and epileptic scenarios~\cite{ullah_influence_2009,barreto_ion_2011}.
\par
The variables in the BC ionic model are: intracellular sodium $\Na$, extracellular potassium $\K$, and intracellular calcium $\Ca$ concentrations. Moreover, three equations are associated with the gating variables $m$ and $h$ that correspond to the activating and inactivating sodium gates, respectively, and $n$ associated to the activating potassium gate. Then, according to the previously introduced notation:
\begin{equation*}
    \boldsymbol{y} = \left[\Ca,\;\K,\;\Na,\;m,\;h,\;n\right]^\top.
\end{equation*}
The Barreto-Cressman ionic model is defined as follows (corresponding to Equation~\eqref{eq:monodomain:ionic}):
\begin{subequations}
\small
\label{eq:barreto_cressman}
    \begin{alignat}{3}
    \label{eq:barreto_cressman:ca}
    \dfrac{\mathrm{d}\Ca}{\mathrm{d}t} & = - \dfrac{\Ca}{\tau_\mathrm{Ca}} 
    - 0.002\,G_\mathrm{Ca}(u - E_\mathrm{Ca})\left(1 + \exp{\left(-\frac{25 + u}{2.5}\right)}\right)^{-1} & & \quad \text{in } (0, T], \\[0pt]
    \dfrac{\mathrm{d}\K}{\mathrm{d} t} & = -\dfrac{1}{\tau} \left(I_\mathrm{diff} + 14I_\mathrm{pump} + I_\mathrm{glia} - 7\gamma I_\mathrm{K}\right) & & \quad \text{in } (0, T], \\[2pt]
    \dfrac{\mathrm{d} \Na}{\mathrm{d}t} & = -\dfrac{1}{\tau} \left(\gamma I_\mathrm{Na} + 3 I_\mathrm{pump}\right) & & \quad \text{in } (0, T], \\[2pt]
    \dfrac{\mathrm{d}g}{\mathrm{d}t} & = \dfrac{3}{\tau_{g}}(g_\infty - g) & & \quad \text{in } (0, T], \quad
    \end{alignat}
\end{subequations}
with $g=\{m,\,h,\,n\}$, $\tau$ is a conversion factor from $\mathrm{s}$ to $\mathrm{ms}$, and $\gamma$ is a unit conversion factor that converts the membrane currents into concentration fluxes. 
\par
The ionic current is characterized by different contributions of sodium, potassium, and chlorine currents, defined as:
\begin{subequations}
\label{eq:current}
    \begin{alignat}{3}
    \label{eq:current:na} 
    I_\mathrm{Na} &= \left(G_\mathrm{NaL} + G_\mathrm{Na}\,m^3 h \right) \left(u - E_\mathrm{Na}\right), \\[5pt]
    \label{eq:current:k} 
    I_\mathrm{K} &= \left(G_\mathrm{K} n^4 + G_\mathrm{AHP} \frac{\Ca}{1 + \Ca} + G_\mathrm{KL} \right) \left(u - E_\mathrm{K}\right), 
    \\[5pt]
    \label{eq:current:cl} 
    I_\mathrm{Cl} &= G_\mathrm{ClL} \left(u - E_\mathrm{Cl}\right).
    \end{alignat}
\end{subequations}
Then, exploiting the definitions in Equations~\eqref{eq:current}, we define the forcing $f(u,\boldsymbol{y})=I_\mathrm{Na}+I_\mathrm{K}+I_\mathrm{Cl}$. We define $(G_{\mathrm{NaL}},\  G_{\mathrm{Na}},\  G_{\mathrm{K}},\  G_{\mathrm{AHP}},\  G_{\mathrm{KL}},\  G_{\mathrm{ClL}})$ as the conductances of the model while the Nerst reversal potentials ($E_\mathrm{Na}$, $E_\mathrm{K}$, $E_\mathrm{Cl}$) are defined as:
\begin{equation*}
\begin{aligned}
E_{\mathrm{Ca}} = & 120\;\mathrm{mV}, && \qquad
E_{\mathrm{K}} = 26.64 \log \left(\frac{\K}{[\mathrm{K}^+]_i} \right) \;\mathrm{mV} ,\\
E_{\mathrm{Na}} = & 26.64 \log \left(\frac{[\mathrm{Na}^+]_0}{\Na} \right) \;\mathrm{mV}, && \qquad 
E_{\mathrm{Cl}} = 26.64 \log \left(\frac{[\mathrm{Cl}^+]_i}{[\mathrm{Cl}^+]_0} \right)\; \mathrm{mV}.  
\end{aligned}
\end{equation*}
The system also includes three differnt types of current that affect the ion concentrations: the one related to the capacity of glial cells to remove excess potassium from the extracellular space ($I_\mathrm{Glia}$), the current that represents the diffusion of potassium ($I_\mathrm{diff}$) and the one related to the sodium–potassium pump ($I_\mathrm{pump}$). We define these currents as:
\begin{equation*}
\begin{aligned}
I_\mathrm{pump} = & \, \rho \left(1 + e^{5.5 - \K}\right)^{-1} \left(1 + \exp\left(\dfrac{25 - \Na}{3}\right) \right)^{-1}, \\[3pt]
I_\mathrm{glia} = & \, G_\mathrm{glia}\left(1 + \exp\left(\dfrac{18 - \K}{2.5}\right)\right)^{-1}, \\[8pt]
I_\mathrm{diff} = & \, \epsilon \left(\K - \Kbath\right),
\end{aligned}
\end{equation*}
where $\rho$, and $\epsilon$ are constants introduced to fit the temporal dynamics of the system, and $\Kbath$ denotes the potassium concentration in the largest nearby reservoir, physiologically associated with the brain vasculature. This parameter plays a crucial role, as an increase in $\Kbath$ can be defined as a trigger for the starting of an epileptic event. For a comprehensive analysis  of the variables, parameters, and assumptions involved in the model, we refer to \cite{cressman_influence_2009}.

\begin{remark}[Calcium dynamics] $\Capl$ contributes to modulate neuronal excitability \cite{steinlein_calcium_2014} and is modeled in Equation~\eqref{eq:barreto_cressman:ca}. We assume that the intracellular calcium concentration is governed by a \textit{leaky-integrator} model that involves both the accumulation of the quantity and its removal over time. The integration captures the intracellular influx or production of calcium ions, which occurs exlusively through voltage-gated calcium channels. In fact, during an epileptic seizure, the marked rise in intracellular calcium is primarily attributable to the influx via these channels. Calcium removal is then mediated through the Plasma Membrane Calcium-ATPase (PMCA) pump.
\end{remark}

\begin{remark}[Calcium-sensitive potassium current]
In the Barreto-Cressman ionic model, the intracellular calcium concentration is exploited specifically to activate a calcium-dependent potassium current, which plays a key role in spike-frequency adaptation in many excitatory neurons. This current flows through potassium channels whose probability of opening increases as intracellular calcium rises. Calcium-sensitive potassium channels contribute to membrane repolarization following an action potential, thereby preparing the neuron for subsequent firing. When intracellular calcium levels rise, these channels open and allow potassium to exit the cell, generating an outward current.
The calcium-gated potassium current are activated When intracellular calcium increases, causing a flow of potassium ions out of the cell. These channels play a key role in controlling the firing rate of the neuron by shaping the time interval between successive action potentials. When their activity is reduced, the resulting afterhyperpolarization is smaller, keeping the membrane potential closer to threshold, making the initiation of another action potential easier and faster.
\end{remark}

\subsection{Impact of amyloid-\textbeta{}: a novel ionic model}
In this section, starting from the BC model, we derive and investigate a novel comprehensive mathematical model for \Abeta{-}mediated multi-pathway intracellular $\Capl$ dynamics. Indeed, in Alzheimer's disease, amyloid-$\beta$ oligomers are known to disturb intracellular calcium regulation  \cite{kuchibhotla_abeta_2008}, leading neuronal excitability \cite{steinlein_calcium_2014}. \Abeta{ } oligomers can affect $\Capl$ exchanges with the extracellular medium and internal $\Capl$ stores \cite{berridge_calcium_2010}.
\par
In \cite{vonbonhorst_impact_2022}, the authors proposed a model involving three $\Capl$ compartments internal to the cell: sub-plasmalemmal space, cytosol, and endoplasmic reticulum (ER). In particular, the sub-plasmalemmal space is a thin fictitious shell just below the plasma membrane that controls its electrical activity. The findings suggest that the $\Capl$ in the sub-plasmalemmal space is governed by the exchange with the extracellular medium, occurring at a much higher rate than with the cytosol and the ER \cite{vonbonhorst_impact_2022}. Consequently, alterations in the activity of calcium transporters induced by \Abeta{ }peptides can significantly affect neuronal excitability. In order to capture these effects, we introduce a novel ionic model, based on the Barreto-Cressman framework, that provides a more accurate representation of calcium $\Capl$ exchange across the membrane. 
\par
\Abeta{-}induced disruptions in intracellular calcium homeostasis can lead either hyperexcitability or hypoexcitability, depending on the mechanisms involved. Computational modeling provides a framework to investigate these effects, particularly the one arising from alterations in PMCA activity, of $\Capl$-permeable membrane pores, the overactivation of voltage-gated $\Capl$ channels, and the inhibition of fast-inactivating $\Kpl$ currents. In physiological conditions, \Abeta{ }accumulation is a slow process that unfolds over months or even decades. Therefore, in our model, we assume a fixed value of \Abeta{ }concentration, representing a specific stage in progression of Alzheimer disease \cite{latulippe_mathematical_2018}. We consider \Abeta{ }concentrations of $0.1\ \muM$, $1\ \muM$, and $10\ \muM$, corresponding to early, intermediate, and advanced stages of the disease \cite{raskatov_what_2019}.
\subsubsection*{Impact of \Abeta{ }peptides on plasma membrane $\Capl$ ATPase (PMCA)}
The PMCA pump is the main mechanism responsible for extracting calcium $\Capl$ from the sub-plasmalemmal space into the extracellular region\cite{vonbonhorst_impact_2022}. This ATP-dependent process is fundamental for maintaining intracellular $\Capl$ homeostasis and for controlling neuronal excitability. 
In Alzheimer disease, however, PMCA activity is impaired due to inhibitory action of \Abeta{ }peptides \cite{berrocal_calmodulin_2012}. As found in \cite{berrocal_altered_2009}, the inhibitory effect of \Abeta{ }follows a inhibition curve described by the equation:
\begin{equation}
\label{eq:pmca_inhibition}
\frac{\varphi^{\mathrm{PMCA}}_{\mathrm{Ca}}}{\varphi^{\mathrm{PMCA}}_{\mathrm{Ca},\mathrm{max}}} = \frac{k_\mathrm{I}}{\Ab + k_\mathrm{I}},
\end{equation}
where $u_{\mathrm{max}}$ and $u$ are ATP hydrolysis rates in absence and presence of \Abeta{ }at a concentration $\Ab$, respectively. The constant $k_\mathrm{I}$ is the apparent dissociation constant for \Abeta{ }peptides, defined as the concentration of \Abeta{ }required to reduce the enzymatic activity of PMCA by $50\%$. The estimated  value of $k_\mathrm{I}$ is $2.312\,\mathrm{\mu M}$ \cite{berrocal_altered_2009}.
\par
The term modelling the PMCA effect on $\Capl$ in Equation~\eqref{eq:barreto_cressman:ca} is $\varphi^{\mathrm{PMCA}}_{\mathrm{Ca},\mathrm{max}} = -\tau_\mathrm{Ca}^{-1}\,\Ca$. According to the findings in \cite{berrocal_altered_2009}, we model the impact of the presence of \Abeta{ }, by defining $\varphi^{\mathrm{PMCA}}_{\mathrm{Ca}} = -(\tau_\mathrm{Ca}+k_\mathrm{PMCA}\Ab)^{-1}\,\Ca$. Substituting these definitions in Equation~\eqref{eq:pmca_inhibition}, we obtain the relation:
\begin{equation}
    k_\mathrm{PMCA} = \frac{\tau_\mathrm{Ca}}{k_\mathrm{I}},
\end{equation}
whose value can be estimated equal to $34.602\,\mathrm{ms}\,\muM^{-1}$.
\subsubsection*{\Abeta{-}induced formation of $\Capl$ permeable PM pores}
One important mechanism by which \Abeta{} peptides are thought to modify $\Capl$ regulation involves the formation of $\Capl$-permeable pores in the membrane, which allow uncontrolled calcium entry into neurons \cite{mirdha_aggregation_2024}. Following the modeling strategy in \cite{vonbonhorst_impact_2022}, we represent the resulting  $\Capl$ influx into the sub-plasmalemmal compartment as:
\begin{equation}
\label{eq:abeta_pores_influx}
J_{A\beta} = J_{A\beta}^\mathrm{max}\left(1 + \exp{\left(\dfrac{u-q_1}{q_2}\right)}\right)^{-1},
\end{equation}
where $J_{A\beta}^\mathrm{max}$ is the maximal rate of $\Capl$ entry through the pores that reflects both the single pore permeability and the number of pores in the PM. Here $q_1 = 30\,\mathrm{mV}$, and $q_2 = 25\,\mathrm{mV}$ are defined as parameters that characterize the voltage-dependence of $\Capl$ influx through the \Abeta{ }pores \cite{vonbonhorst_impact_2022}.
\par
 At $u \approx q_1$, the $\Capl$ influx reaches approximately half of its maximal value $ J_{A\beta}^\mathrm{max}$. We describe the dependence of the influx rate $V_{A\beta}$ on the concentration of \Abeta{}using a Hill-type relation of the form:
\begin{equation}
\label{eq:pores_rate}
 J_{A\beta}^\mathrm{max} = J_{A\beta}^\mathrm{asy}\frac{\Ab}{k_\mathrm{D} + \Ab},
\end{equation}
where $J_{A\beta}^\mathrm{asy}$ is the asymptotic maximum rate of the channels, which is achieved when the concentration of \Abeta{ } is sufficiently high, and we fix it equal to $J_{A\beta}^\mathrm{asy}=10\;\muM\,\mathrm{ms}^{-1}$ \cite{vonbonhorst_impact_2022}. Moreover, we set $k_\mathrm{D} = 10\,\muM$. With this formulation, an increase in \Abeta{ }concentration leads to a higher maxima rate of calcium entry through the pores, thereby elevating the sub-plasmalemmal calcium level.
\subsubsection*{Overactivation of voltage-gated $\Capl$ channels induced by \Abeta{ }peptides}
Disruptions in intracellular calcium homeostasis can also arise from abnormal calcium $\Capl$ influx through plasma membrane channels, including voltage-gated $\Capl$ channels (VGCCs). Here, L-type $\Capl$ channels are strongly implicated in the neurotoxic effects of \Abeta{ }peptides \cite{kim_effects_2011}. 
Experimental studies in rat models have shown that, in the presence of oligomeric \Abeta{-42},these channels become over-activated and open at lower membrane potentials than under normal conditions \cite{ishii_amyloid_2019}. This results in an enhanced calcium influx, causing intracellular calcium to rise more rapidly and reach higher levels than under physiological conditions. To account for the influence of \Abeta{ }on L-type VGCCs, we modify the activation dynamics of these channels in the Barreto-Cressman ionic model. The resulting forcing term is:
\begin{equation}
\small
\begin{split}
I_\mathrm{VGCC} = &
    -  0.002 G_\mathrm{Ca}(u - E_\mathrm{Ca})\left(1 + \exp\left(-\frac{1}{2.5}\left(25 + u + \frac{u_{A\beta}^\mathrm{max}\Ab^\alpha}{k_\mathrm{VGCC}^\alpha+\Ab^\alpha}\right)\right)\right)^{-1} \\
    = & -  0.002 G_\mathrm{Ca}(u - E_\mathrm{Ca})\left(1 + \exp\left(-\frac{25 + u + u_{A\beta}}{2.5}\right)\right)^{-1},
\end{split}
\end{equation}
here $\alpha = 0.5$ and $u_{A\beta}^\mathrm{max}$ denotes the maximum \Abeta{-}induced shift in the activation voltage of the channel. Experimental evidence shows that an \Abeta{-42} concentration of $0.1 \ \muM$ induces a shift of approximately $-20 \mathrm{mV}$ in the activation threshold of L-type VGCCs \cite{ishii_amyloid_2019}. 
In the standard model, the half-activation voltage occurs at $u=-25\ \mathrm{mV}$, meaning the channel is open with $50\%$ probability at this potential.
With $u_{A\beta}^\mathrm{max} = 25 \ \mathrm{mV}$, the channels has a $50\%$ chance of opening at $u= - 50\ \mathrm{mV}$ allowing earlier $\Capl$ influx. By using these literature data, we fix the parameter $k_{\mathrm{VGCC}}=4.44\times10^{-2} \muM$.
\subsubsection*{\Abeta{-}mediated block of the fast-inactivating potassium current}
\Abeta{ }peptides have been shown to impair the fast-inactivating, voltage-gated potassium current, reducing potassium efflux and thereby prolonging depolarization and enhancing 
calcium entry \cite{good_effect_1996}. 
Evidence from rat studies indicates that \Abeta{ }binds to the closed state of these channels and prevents them from opening, without altering their intrinsic gating kinetics or voltage dependence \cite{good_effect_1996}. To incorporate this effect, we express the maximal conducance of the fast-inactivating potassum channels as the conducance in the absence of \Abeta{,} scaled by the fraction of channels that remain available to open. Following the idea in \cite{good_effect_1996}, the new expression in the $I_\mathrm{K}$ definition is:
\begin{equation}
    G_\mathrm{K} \left(1-\dfrac{\Ab}{k_\mathrm{Ca,K}+\Ab}\right) n^4,
\end{equation}
where $k_\mathrm{Ca,K}$ is the inhibition constant with an estimated value of $10\,\muM$ \cite{good_beta-amyloid_1996}.
\subsubsection*{\Abeta{-}mediated decrease the activity of $\Capl$-sensitive potassium channels}
The activity of the $\Capl$-sensitive potassium channels has also been shown to be modified by the presence of \Abeta{ }peptides, making them less active \cite{yamamoto_suppression_2011}. When those channels' conductance decreased, their after-hyperpolarizing effect is reduced, and they become less active at 
allowing potassium ions to flow through \cite{yamamoto_amyloid_2021}. As a consequence, it is easier to trigger repetitive firing. Indeed, \Abeta{-42} increases the calcium influx thanks to the 
suppression of BK channels, which leads to 
prolonged spike duration and thus increased Ca2+ entry. Indeed, BK channels are crucial for quickly depolarizing the membrane after depolarization. Less activity of those channels slows the depolarization, leading to a border spike. The reduction of current produced by those channels is estimated around $8.77\%$ with $\Ab=1\,\muM$ \cite{jhanandas_cellular_2001}, and around $20.80\%$ with $\Ab=5\,\muM$ \cite{zhang_effects_2014}. For this reason, we propose a scaling of the BK currents associated with $G_\mathrm{AHP}$ by a factor defined as
\begin{equation*}
    S_{A\beta} = \dfrac{a_\mathrm{BK}}{\Ab+b_\mathrm{BK}} + c_\mathrm{BK}\,e^{-d_\mathrm{BK}\,\Ab}
\end{equation*}
where $a_\mathrm{BK} = 4.498\times 10^{-1} \muM$, $b_\mathrm{BK} = 1.9295\,\muM$, $c_\mathrm{BK} = 0.7669$, and $d_\mathrm{BK} = 10.70\,\muM^{-1}$.
\subsubsection{The \Abeta{-}dependent neuronal model}

\begin{figure}[h]
    \centering
\includegraphics[width=\textwidth]{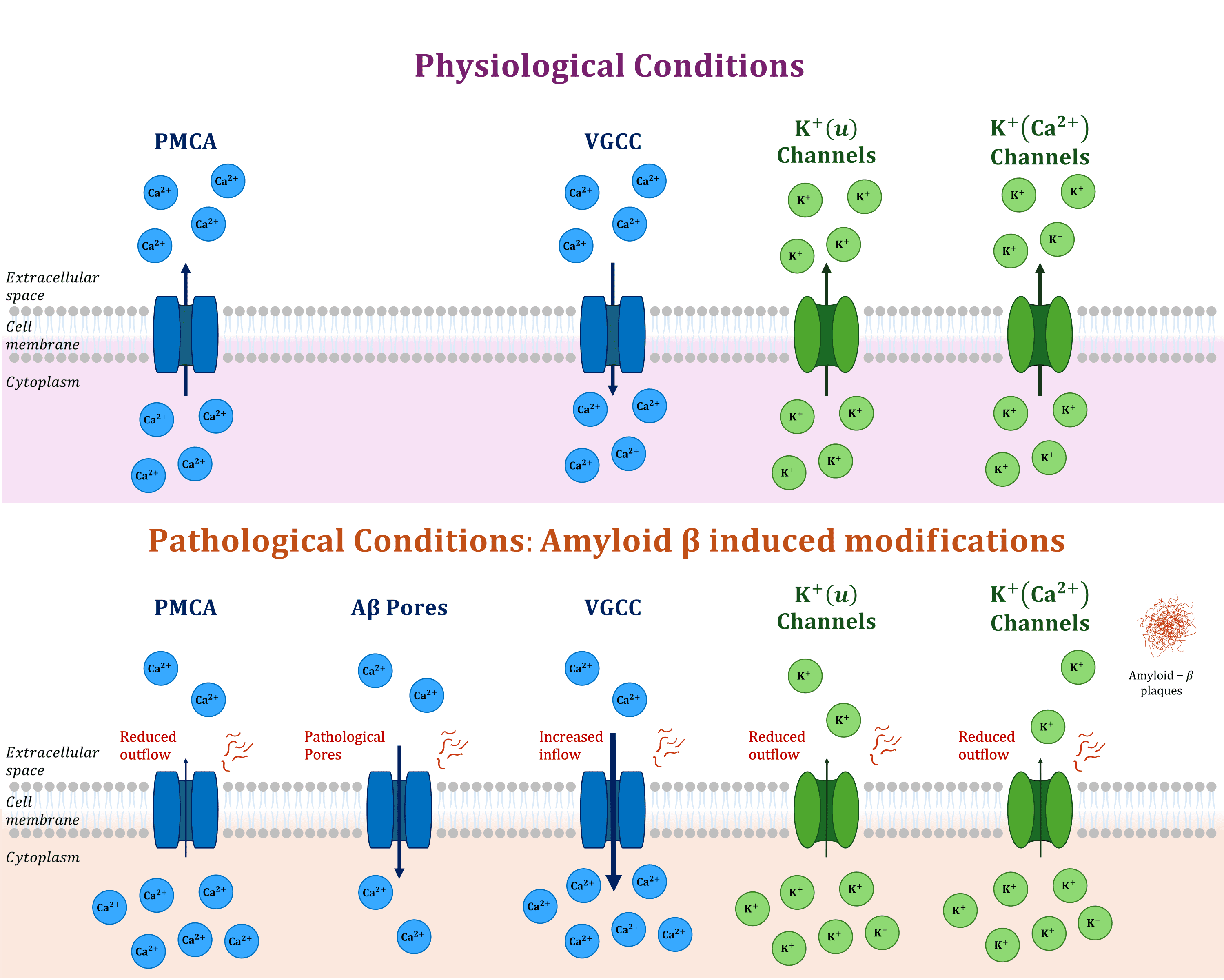}
    \caption{Schematic model of the pathological modifications induced by the \Abeta{ }peptides.}
    \label{fig:criticalfluxes}
\end{figure}

In Figure \ref{fig:criticalfluxes}, we summarize the critical fluxes induced by \Abeta{ }peptides and modeled in this section. Finally, our model, depending on the $\Ab$ concentration, consists of the following system of ODEs (corresponding to Equation~\eqref{eq:monodomain:ionic}):
\begin{subequations}
\small
\label{eq:bc_abeta}
    \begin{alignat}{4}
\label{eq:bc_abeta:ca}
\nonumber \dfrac{\mathrm{d}\Ca}{\mathrm{d}t}
= {} & - \dfrac{\Ca}{\tau_\mathrm{Ca}+k_\mathrm{PMCA}\Ab}  + J_{A\beta} \\
     & - 0.002\,G_\mathrm{Ca}(u - E_\mathrm{Ca})\left(1 + \exp\left(-\frac{25 + u + u_{A\beta}}{2.5}\right)\right)^{-1}
&&\quad \text{in } (0, T], \\[0pt]
    \label{eq:bc_abeta:k}
    \dfrac{\mathrm{d}\K}{\mathrm{d} t} & = -\dfrac{1}{\tau} \left(I_\mathrm{diff} + 14I_\mathrm{pump} + I_\mathrm{glia} - 7\gamma I_\mathrm{K}\right) & & \quad \text{in } (0, T], \\[2pt]
    \label{eq:bc_abeta:na}
    \dfrac{\mathrm{d} \Na}{\mathrm{d}t} & = -\dfrac{1}{\tau} \left(\gamma I_\mathrm{Na} + 3 I_\mathrm{pump}\right) & & \quad \text{in } (0, T], \\[2pt]
    \label{eq:bc_abeta:gating}
    \dfrac{\mathrm{d}g}{\mathrm{d}t} & = \dfrac{3}{\tau_{g}}(g_\infty - g) & & \quad \text{in } (0, T],
    \end{alignat}
\end{subequations}
with $g=m,\,h,\,n$. Moreover, in our model $f(u,\boldsymbol{y})=I_\mathrm{Na}+I_\mathrm{K}+I_\mathrm{Cl}+ \gamma^{-1} J_{A\beta}$, where $\gamma$ is used as a scaling parameter (see \cite{cressman_influence_2009}). Moreover, the modified ionic currents read:
\begin{subequations}
\small
\label{eq:new_current}
    \begin{alignat}{3}
    \label{eq:new_current:na} 
    I_\mathrm{Na} &= \left(G_\mathrm{NaL} + G_\mathrm{Na}\,m^3 h \right) \left(u - E_\mathrm{Na}\right), \\[5pt]
    \label{eq:new_current:k} 
    I_\mathrm{K} &= \left(G_\mathrm{K} \left(1-\dfrac{\Ab}{k_\mathrm{Ca,K}+\Ab}\right) n^4 + S_{A\beta}G_\mathrm{AHP} \frac{\Ca}{1 + \Ca} + G_\mathrm{KL} \right) \left(u - E_\mathrm{K}\right), 
    \\[5pt]
    \label{eq:new_current:cl} 
    I_\mathrm{Cl} &= G_\mathrm{ClL} \left(u - E_\mathrm{Cl}\right).
    \end{alignat}
\end{subequations}

All the novel introduced parameters, together with their biological interpretation, values, and units, are summarized in Table~\ref{tab:model_parameters}.

\begingroup
\small
\setlength{\tabcolsep}{4pt}
\renewcommand{\arraystretch}{1.15}
\setlength{\LTpost}{0pt}

\begin{longtable}{
    @{}
    >{\raggedright\arraybackslash}p{0.16\textwidth}
    >{\raggedright\arraybackslash}p{0.46\textwidth}
    >{\raggedright\arraybackslash}p{0.10\textwidth}
    >{\raggedright\arraybackslash}p{0.10\textwidth}
    >{\raggedright\arraybackslash}p{0.10\textwidth}
    @{}
}

\toprule
\textbf{Symbol} &
\textbf{Biological interpretation} &
\textbf{Value} &
\textbf{Unit} &
\textbf{Reference} \\
\midrule
\endfirsthead

\multicolumn{5}{c}{\tablename~\thetable{} -- continued}\\
\toprule
\textbf{Symbol} &
\textbf{Biological interpretation} &
\textbf{Value} &
\textbf{Unit} &
\textbf{Reference} \\
\midrule
\endhead

\bottomrule
\endfoot

\endlastfoot

\arrayrulecolor{gray!45}
\specialrule{0.3pt}{6pt}{4pt}
\arrayrulecolor{black}

\multicolumn{5}{@{}l}{\textbf{\Abeta{-}dependent PMCA inhibition}}\\*

\arrayrulecolor{gray!45}
\specialrule{0.3pt}{6pt}{4pt}
\arrayrulecolor{black}

$k_{\mathrm{I}}$ &
\Abeta{} concentration producing
half-maximal inhibition of PMCA activity &
$2.312$ & $\mu\mathrm{M}$ &
\cite{berrocal_altered_2009} \\

$k_{\mathrm{PMCA}}$ &
Increase in the calcium-removal time constant
per unit \Abeta{} concentration &
$34.602$ & $\mathrm{ms}\,\mu\mathrm{M}^{-1}$ &
\cite{berrocal_altered_2009} \\

\arrayrulecolor{gray!45}
\specialrule{0.3pt}{6pt}{4pt}
\arrayrulecolor{black}

\multicolumn{5}{@{}l}{\textbf{Amyloid-$\beta$-induced calcium-permeable pores}}\\*

\arrayrulecolor{gray!45}
\specialrule{0.3pt}{6pt}{4pt}
\arrayrulecolor{black}

$J_{A\beta}^{\mathrm{asy}}$ &
Asymptotic maximal calcium influx through
amyloid-$\beta$-induced pores &
$10$ & $\mu\mathrm{M}\,\mathrm{ms}^{-1}$ &
\cite{vonbonhorst_impact_2022} \\

$k_{\mathrm{D}}$ &
Amyloid-$\beta$ concentration at which
$J_{A\beta}^{\max}=J_{A\beta}^{\mathrm{asy}}/2$ &
$10$ & $\mu\mathrm{M}$ &
\cite{vonbonhorst_impact_2022} \\

$q_1$ &
Membrane potential at which pore-mediated
calcium influx equals $J_{A\beta}^{\max}/2$ &
$30$ & $\mathrm{mV}$ &
\cite{vonbonhorst_impact_2022} \\

$q_2$ &
Voltage scale controlling the steepness
of pore-mediated calcium influx &
$25$ & $\mathrm{mV}$ &
\cite{vonbonhorst_impact_2022} \\

\arrayrulecolor{gray!45}
\specialrule{0.3pt}{6pt}{4pt}
\arrayrulecolor{black}

\multicolumn{5}{@{}l}{\textbf{Amyloid-$\beta$-dependent calcium-channel activation}}\\*

\arrayrulecolor{gray!45}
\specialrule{0.3pt}{6pt}{4pt}
\arrayrulecolor{black}

$u_{A\beta}^{\max}$ &
Maximal magnitude of the shift towards
lower calcium-channel activation potentials &
$25$ & $\mathrm{mV}$ &
\cite{ishii_amyloid_2019} \\

$k_{\mathrm{VGCC}}$ &
Amyloid-$\beta$ concentration producing
half of the maximal activation-voltage shift &
$4.44\times10^{-2}$ & $\mu\mathrm{M}$ &
\cite{ishii_amyloid_2019} \\

$\alpha$ &
Hill exponent governing the dependence
of the voltage shift on amyloid-$\beta$ &
$0.5$ & -- &
\cite{ishii_amyloid_2019} \\

\arrayrulecolor{gray!45}
\specialrule{0.3pt}{6pt}{4pt}
\arrayrulecolor{black}

\multicolumn{5}{@{}l}{\textbf{Amyloid-$\beta$-dependent potassium-channel inhibition}}\\*

\arrayrulecolor{gray!45}
\specialrule{0.3pt}{6pt}{4pt}
\arrayrulecolor{black}

$k_{\mathrm{Ca,K}}$ &
Amyloid-$\beta$ concentration producing
half-maximal inhibition of the $G_{\mathrm{K}}$ contribution &
$10$ & $\mu\mathrm{M}$ &
\cite{good_beta-amyloid_1996} \\

$a_{\mathrm{BK}}$ &
Numerator coefficient of the rational contribution
to the calcium-sensitive potassium-current scaling &
$4.498\times10^{-1}$ & $\mu\mathrm{M}$ &
\cite{jhanandas_cellular_2001,zhang_effects_2014} \\

$b_{\mathrm{BK}}$ &
Concentration scale of the rational contribution
to the calcium-sensitive potassium-current scaling &
$1.9295$ & $\mu\mathrm{M}$ &
\cite{jhanandas_cellular_2001,zhang_effects_2014} \\

$c_{\mathrm{BK}}$ &
Amplitude of the exponential contribution
to the calcium-sensitive potassium-current scaling &
$0.7669$ & -- &
\cite{jhanandas_cellular_2001,zhang_effects_2014} \\

$d_{\mathrm{BK}}$ &
Concentration-dependent decay coefficient
of the exponential contribution &
$10.70$ & $\mu\mathrm{M}^{-1}$ &
\cite{jhanandas_cellular_2001,zhang_effects_2014} \\*

\bottomrule

\caption{New parameters of the \Abeta{-}dependent ionic model.
The coefficients $a_{\mathrm{BK}}$, $b_{\mathrm{BK}}$,
$c_{\mathrm{BK}}$, and $d_{\mathrm{BK}}$ are estimated
in this work based on the experimental data
reported in the corresponding references.}
\label{tab:model_parameters}\\

\end{longtable}
\endgroup

\section{PolyDG formulation}
\label{sec:3}
We now present the PolyDG semi-discrete formulation of the problem described in Equation \eqref{eq:monodomain}, which is derived in detail in \cite{leimer_saglio_high-order_2024}. Let $\mathcal{T}_h$ represent a polytopal mesh partition of the domain $\Omega$, consisting of disjoint elements $K$. For each element $K$, we define its diameter as $h_K$ and set $h = \max_{K \in \mathcal{F}_h} h_K < 1$. The interfaces are defined as the intersections of the $(d-1)$-dimensional facets of neighboring elements. We denote by $\mathcal{F}_h^I$ the union of all interior faces contained within $\Omega$ and by $\mathcal{F}_h^N$ those lying on the boundary $\partial \Omega$. In the following, we assume that the underlying grid is polytopic regular in the sense of \cite{cangiani_hp-version_2014,cangiani_hp-version_2017}.
\par
We define $\mathbb{P}^{p_K}(K)$ as the space of polynomials of degree $p_K\geq1$ over the element $K$ and the discontinuous finite element space as: 
\begin{equation*}
V_h^{\mathrm{DG}} = \{v_h \in L^2(\Omega) : v_h|_{K} \in \mathbb{P}^{p_K}(K) \quad \forall \: K \in \mathcal{T}_h\},
\end{equation*}
Let $F \in \mathcal{F}_h^I$ be the face shared by the elements $K^{\pm}$, and let $\boldsymbol{n}^{\pm}$ denote the normal unit vectors pointing outward to $K^\pm$, respectively. For a regular enough scalar-valued function $v$ and a vector-valued function $\boldsymbol{q}$, the trace operators are defined as follows \cite{arnold_unified_2002}:
\begin{equation*}
    \begin{aligned}
      &\averagel v \averager = \frac{1}{2} (v^+ + v^-) , \quad & \jumpl v  \jumpr& = v^+ \boldsymbol{n}^+ + v^- \boldsymbol{n}^-, \quad &\text{on }F \in \mathcal{F}_h^I,&\\
      &\averagel \boldsymbol{q} \averager =  \frac{1}{2} (\boldsymbol{q}^+ + \boldsymbol{q}^-) , \quad & \jumpl \boldsymbol{q} \jumpr & = \boldsymbol{q}^+\cdot \boldsymbol{n}^+ + \boldsymbol{q}^-\cdot \boldsymbol{n}^-,  \quad &\text{on }F \in \mathcal{F}_h^{I},&\\
    \end{aligned}
\end{equation*}
where the superscripts $\pm$ indicate the traces of these functions on $F$ taken in the interiors of $K^{\pm}$, respectively. The definition of the penalization parameter reads as follows $\eta : \mathcal{F}_h^I \cup \mathcal{F}_h^D \rightarrow \mathbb{R}_+$:
\begin{equation}
    \eta = \eta(\boldsymbol{p},h,\boldsymbol{\Sigma}) = \eta_0 
    \begin{cases}
        \{\boldsymbol{\Sigma}_K\}_A \dfrac{\{p^2_K\}_A}{\{h_K\}_H} & \text{on } F \in \mathcal{F}_h^I, \\
        \boldsymbol{\Sigma}_K \dfrac{p^2_K}{h_K} & \text{on } F \in \mathcal{F}_h^D, \\
    \end{cases}
    \label{eq:eta}
\end{equation}
which depends explicitly on both the local degrees and the mesh size. This allows us to exploit the $p$-adaptive algorithm for travelling wavefront scenarios as described in \cite{leimer_saglio_p-adaptive_2025}.
In Equation \eqref{eq:eta} we set $\boldsymbol{\Sigma}_K = \|\sqrt{\boldsymbol{\Sigma}}|_K\|^2_{L^2(K)}$ and we consider both the harmonic average operator $\{\cdot\}_H$, and the arithmetic average operator $\{\cdot\}_A$ on $F \in \mathcal{F}_h^I$. 
The parameter $\eta_0$ is chosen large enough to ensure stability. This setting enable us to introduce the following bilinear form $\mathcal{A}(\cdot,\cdot): V_h^{\mathrm{DG}}\times V_h^{\mathrm{DG}} \rightarrow \mathbb{R}$:
\begin{equation}
    \begin{aligned}
       \mathcal{A}(u,v) =& \int_{\Omega} \mathbf{\Sigma}\nabla_h u \cdot  \nabla_h v \;dx+ \sum_{F \in \mathcal{F}_h^I} \int_F (\eta  \jumpl u \jumpr \cdot  \jumpl v  \jumpr) \\- &\sum_{F \in \mathcal{F}_h^I} \int_F( \averagel \boldsymbol{\Sigma} \nabla u \averager \cdot  \jumpl v\jumpr   - \jumpl u \jumpr \cdot \averagel \boldsymbol{\Sigma} \nabla v\averager) d\sigma \quad \forall \: u,v \in V^{\text{DG}}_h ,
    \end{aligned}
    \label{eq:coer}
\end{equation}
where $\nabla_h$ is the element-wise gradient. The semi-discrete formulation of problem in Equation \eqref{eq:monodomain} reads:
\par
For any $t \in (0,T]$, find $(u_h(t),\boldsymbol{y}_h(t)) \in V^{\mathrm{DG}}_h \times \mathbf{V}^{\mathrm{DG}}_h$ such that:
\begin{equation}
\label{eq:semi-discrete_barreto}
    \begin{aligned}
         \left(\chi_m C_m \frac{\partial u_h(t)}{\partial t},v_h\right)_{\Omega} +  \mathcal{A}(u_h(t),v_h) \\ \qquad \qquad + \left(\chi_m f(u_h(t),\boldsymbol{y}_h(t)),v_h\right)_{\Omega} &=   \,(I_h^{\mathrm{ext}},v_h)_{\Omega}   &\forall \,  v_h \in V_h^{\mathrm{DG}}, \\
         \left(\frac{\partial\boldsymbol{y}_h(t)}{\partial t},\boldsymbol{w}_h\right)_{\Omega} + \left( \boldsymbol{m}(u_h(t),\boldsymbol{y}_h(t)),\boldsymbol{w}_h\right)_{\Omega} &=  \,0  &\forall \,  \boldsymbol{w}_h \in \mathbf{V}_h^{\mathrm{DG}}, \\
        u_h(0) &= u_h^0,\quad \boldsymbol{y}_h(0) =  \,\boldsymbol{y}_h^0  &\text{in } \Omega.
    \end{aligned}
\end{equation}

As a consequence of our assumptions, we denote the dimension of the discrete space as $N_h(\boldsymbol{p})$, to make explicit its dependence on the vector of local polynomial distribution.
Let $N_h(\boldsymbol{p})$ be the dimension of $V_h^\mathrm{DG}$ and let $(\varphi_j)^{N_h(\boldsymbol{p})}_{j=0}$ be a suitable basis for $V_{h}^\mathrm{DG}$, then $u_h(t) = \sum_{j=0}^{N_h(\boldsymbol{p})} U_j(t)\varphi_j$ and $y_l(t) = \sum_{j=0}^{N_h(\boldsymbol{p})} Y^l_j(t)\varphi_j$ for all $l=1,...,n$. We denote $\mathbf{U} \in \mathbb{R}^{N_h(\boldsymbol{p})}$, $\mathbf{Y}_l \in \mathbb{R}^{N_h(\boldsymbol{p})}$ for all $l=1,...,n$ and $\mathbf{Y} = [\mathbf{Y}_1,...,\mathbf{Y}_n]^\top$. We define the matrices:
\begin{equation}
\small
    \begin{aligned}
      [\mathbf{M}]_{ij} &= (\varphi_i,\varphi_j)_{\Omega}, \; &\text{(Mass matrix),}& \quad  i,j = 1,...,N_h(\boldsymbol{p}) \\
      [\mathbf{F}]_{j} &= (I^\text{ext},\varphi_j)_{\Omega},  \; &\text{(Forcing term),}& \quad j = 1,...,N_h(\boldsymbol{p})\\
      [\mathbf{I}(u,\boldsymbol{y})]_{j} &= (f(u,\boldsymbol{y}),\varphi_j)_{\Omega},  \; &\text{(Non-linear ionic forcing term),}& \quad  j = 1,...,N_h(\boldsymbol{p})\\
      [\mathbf{G}_l(u,\boldsymbol{y})]_{j} &= (\boldsymbol{m}_l(u,\boldsymbol{y}),\boldsymbol{\varphi}_j)_{\Omega},  \; &\text{(Dynamics of the ionic model),}& \quad  j = 1,...,N_h(\boldsymbol{p}), \\ & & & \quad l=1,...,n  \\ 
      [\mathbf{A}]_{ij} &= \mathcal{A}(\varphi_i,\varphi_j) &\text{(Stiffness matrix),}& \quad  i,j = 1,...,N_h(\boldsymbol{p}).
      \label{eq::matrixFull}
    \end{aligned}
\end{equation}
Finally, we introduce the fully-discrete formulation.
We partition the interval \([0, T]\) into \(N\) sub-intervals \((t^{(k)}, t^{(k+1)}]\), each of length \(\Delta t\), such that \(t^{(k)} = k\Delta t\) for \(k = 0, \dots, N-1\). Concerning the time discretization, we adopt the second-order Crank-Nicolson scheme for the linear part, with the ion term discretized with a second-order explicit extrapolation. Given the initial conditions \(\mathbf{U}_0\) and \(\mathbf{Y}_0\), the discrete scheme is: find $\mathbf{U}^{(k+1)} = \mathbf{U}(t^{(k+1)})$ and  $\mathbf{Y}^{(k+1)} = \mathbf{Y}(t^{(k+1)})$ for $k=0,...,N-1$, such that
\begin{equation*}
\label{eq:Full_discrete_complete}
    \begin{aligned}
      \left(\chi_m C_m \mathbf{M} + \frac{\Delta t}{2}\mathbf{A} \right) \mathbf{U}^{(k+1)}  = & \,  \left(\chi_m C_m \mathbf{M} - \frac{\Delta t}{2}\mathbf{A} \right) \mathbf{U}^{(k)}  - \frac{\mathrm{\chi_m} \Delta t}{2} (3\mathbf{I}^{(k)}-\mathbf{I}^{(k-1)}) \\&+ \frac{\Delta t}{2}( \mathbf{F}^{(k+1)}+\mathbf{F}^{(k)}) , \\
      \mathbf{Y}^{(k+1)} = & \, \mathbf{Y}^{(k)} - \Delta t \mathbf{G}^{(k)},\\
      (\mathbf{U}^{(0)},\mathbf{Y}^{(0)}) = &\, (\mathbf{U}_0,\mathbf{Y}_0).
    \end{aligned}
\end{equation*}

\section{Sensitivity analysis of the novel 0D ionic model with respect to amyloid-\textbeta{ }concentration}
\label{sec:0d}
In this section, we investigate the effects of different levels of A\textbeta{ }peptides on calcium dynamics and neuronal excitability using our mathematical model.
\par
\begin{figure}[ht]
     \centering
     \includegraphics[width=\textwidth]{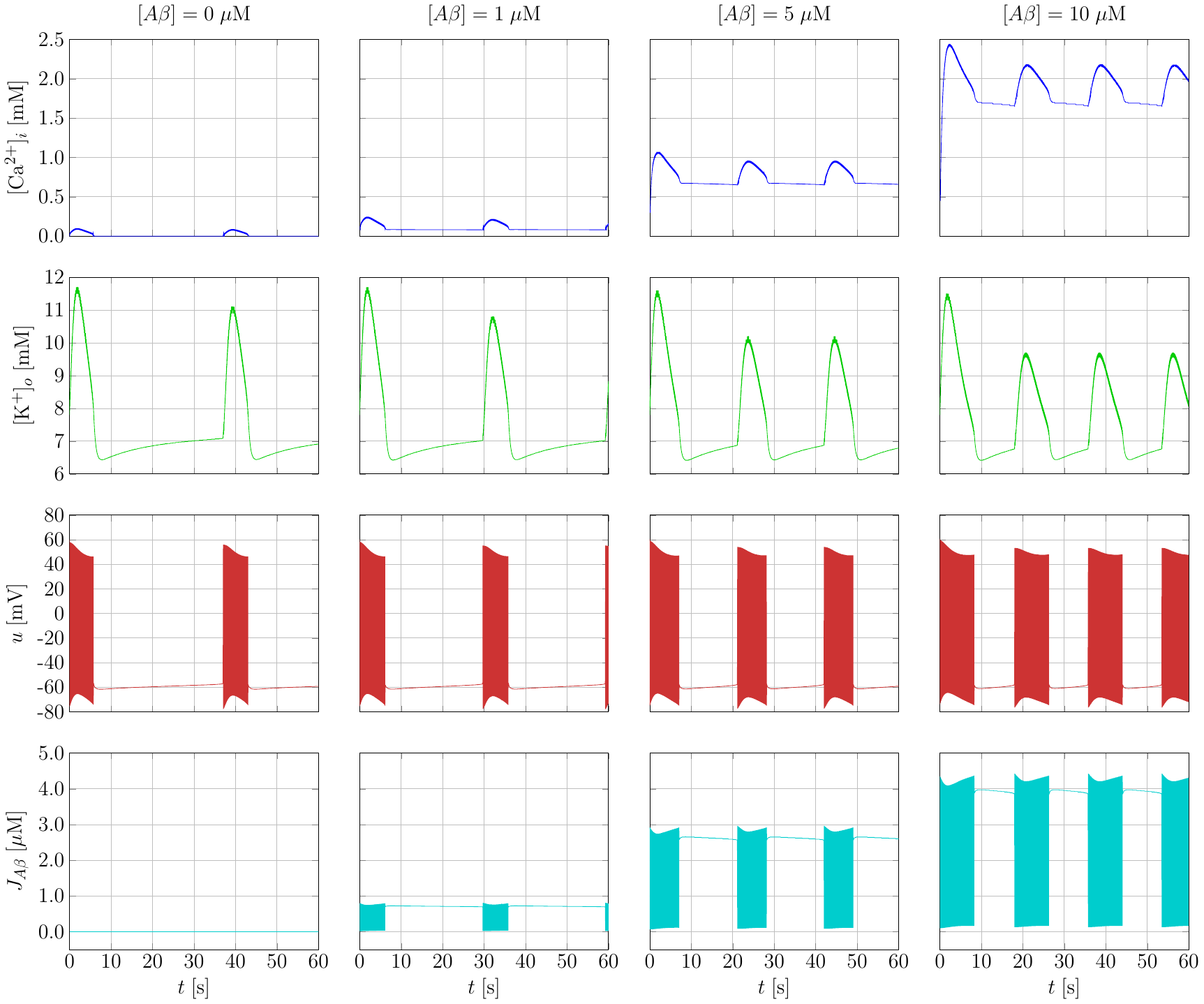}
     \caption{Effects of A\textbeta{ }accumulation with membrane potential with $\mathrm{K}_\mathrm{bath}=8\ \mathrm{mM}$. Evolution over $60 \text{s}$ of $\Ca$ (first row), $\K$ (second row), $u$ (third row) and $J_{A\beta}$ (fourth row) for different values of $\Ab$.
    }
     \label{fig:phase_plots}
\end{figure}
First of all, we perform a simulation with a time interval $(0\,\mathrm{s}, 60\,\mathrm{s})$ for different values of \Abeta{ }concentrations: $\Ab = 0, 1, 5, 10 \, \mu\mathrm{M}$. In Figure \ref{fig:phase_plots}, we report the results of the simulations, showing that the multiple \Abeta-affected pathways introduced in the model contribute to $\Capl$ dysregulation and the subsequent neuronal hyper-excitability.
\par
The first column of Figure \ref{fig:phase_plots} reports a spiking epileptic activity in the absence of \Abeta. For $\Ab=1\,\mu\mathrm{M}$ (second column of Figure \ref{fig:phase_plots}), a slight increase in $\Ca$ concentration is observed. This induces bursts in the membrane potential characterized by a higher frequency and longer lasting, which reflects an altered neuronal excitability. Increasing the concentration until $\Ab = 10\,\mu\mathrm{M}$  (fourth column of Figure \ref{fig:phase_plots}), burst acceleration and prolongation becomes more evident. We can observe a significant increase in $\Ca$, also in the non-spiking times, and a decrease in the concentration of $\K$ in its peak values.  
\par
\begin{figure}[ht]
\centering
\includegraphics[width=\textwidth]{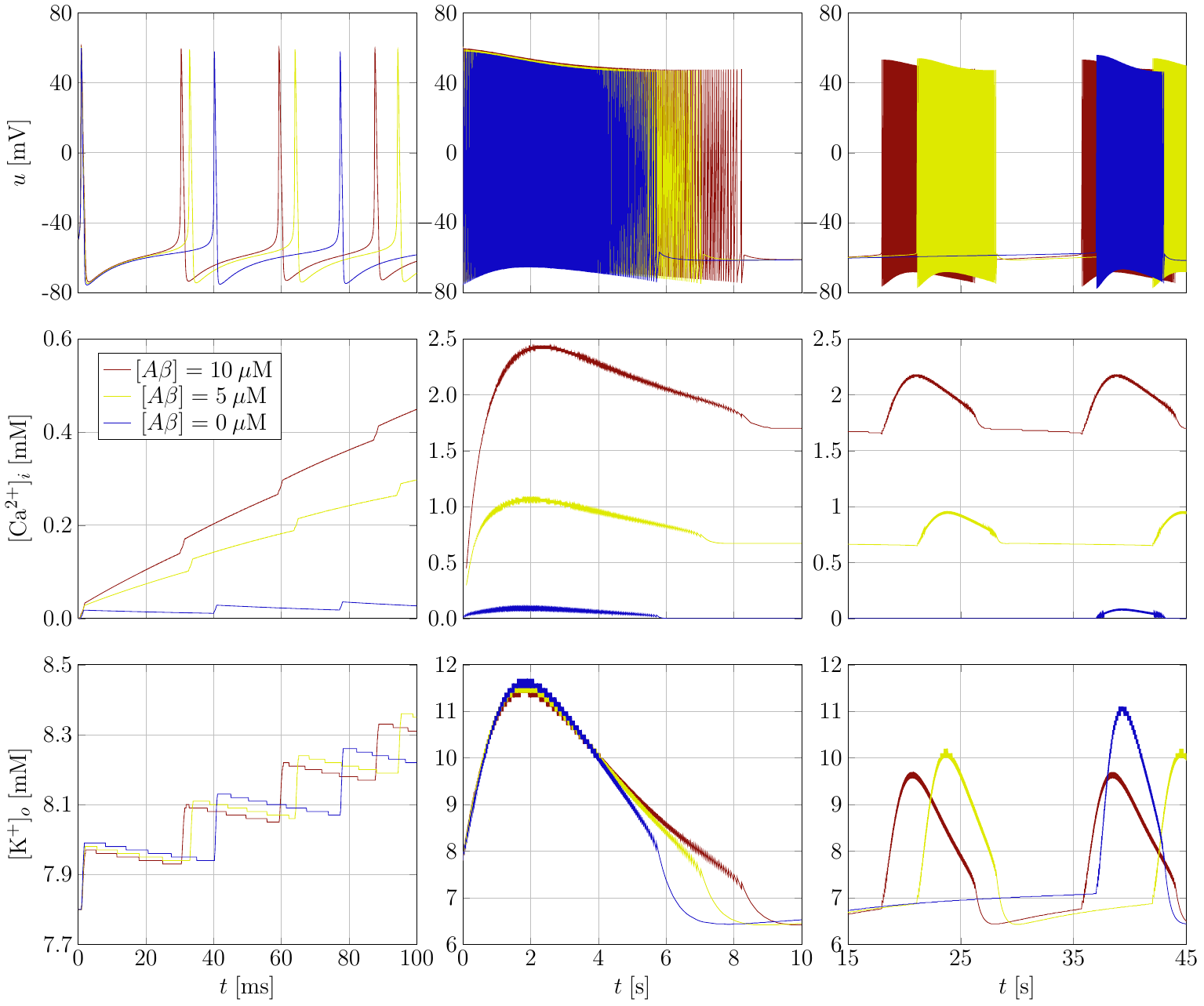}
 \caption{Effects of \Abeta{} accumulation with membrane potential with $\mathrm{K}_\mathrm{bath}=8\ \mathrm{mM}$. Evolution of $\Ca$ (first row), $\K$ (second row), $u$ (third row) for different values of $\Ab$ and on different time scales $(0\,\mathrm{ms},100\,\mathrm{ms})$ (left), $(0\,\mathrm{s},10\,\mathrm{s})$ (center), $(15\,\mathrm{s},45\,\mathrm{s})$ (right).
    }
\label{fig:sensitivity_abeta_2D}
\end{figure}

The results show the burst progression as the concentration of \Abeta{ }increases, with consequent increase in the induced currents $J_{A\beta}$. To better highlight the sensitivity of the quantities with respec to \Abeta{ } concentrations, we report the simulations over different time periods for $\Ab = 0, 5, 10 \, \mu\mathrm{M}$ in Figure \ref{fig:sensitivity_abeta_2D}. In the left column, we show the concentration dynamics in the first $10\,\mathrm{ms}$. In particular, we can notice that $\Ca$ concentration increases faster if there is a high concentration of \Abeta{.} This fact induces a faster spiking of the action potential. In the second and third columns of Figure \ref{fig:sensitivity_abeta_2D}, we report the solution in the intervals $(0\,\mathrm{s},10\,\mathrm{s})$ and $(15\,\mathrm{s},45\,\mathrm{s})$, respectively. On this longer scale, the intracellular $\Ca$ concentration increases sensibly during the initial seconds and then decreases, but with persistent oscillations at high levels, resulting in a pathological and abnormal concentration of calcium in the region affected by $\Ab$ concentration. This behavior reflects a pronounced disruption of calcium homeostasis, in which impaired clearance mechanisms fail to restore intracellular calcium to physiological levels. Such dysregulation has been widely documented in the medical literature \cite{berridge_calcium_2010,kuchibhotla_abeta_2008}. A similar alteration is observed in the membrane potential, which no longer displays isolated bursts but instead remains in a state of persistent, high-frequency spiking. 

The oligomers of \Abeta{ }peptides induce a spike activity enhancement at the neuronal level, consistent with the medical literature findings \cite{yang_alzheimer_2022,romoli_amyloid_2021}. Finally, we can observe a sensible reduction of the peaks of $\K$ concentration, with the increase of \Abeta{,} as well as a slower discharge of the cells after each burst \cite{yamamoto_suppression_2011,yamamoto_amyloid_2021}.
\par

\begin{figure}[t]
    \begin{subfigure}[t]{0.50\textwidth}
     \includegraphics[width=\textwidth]{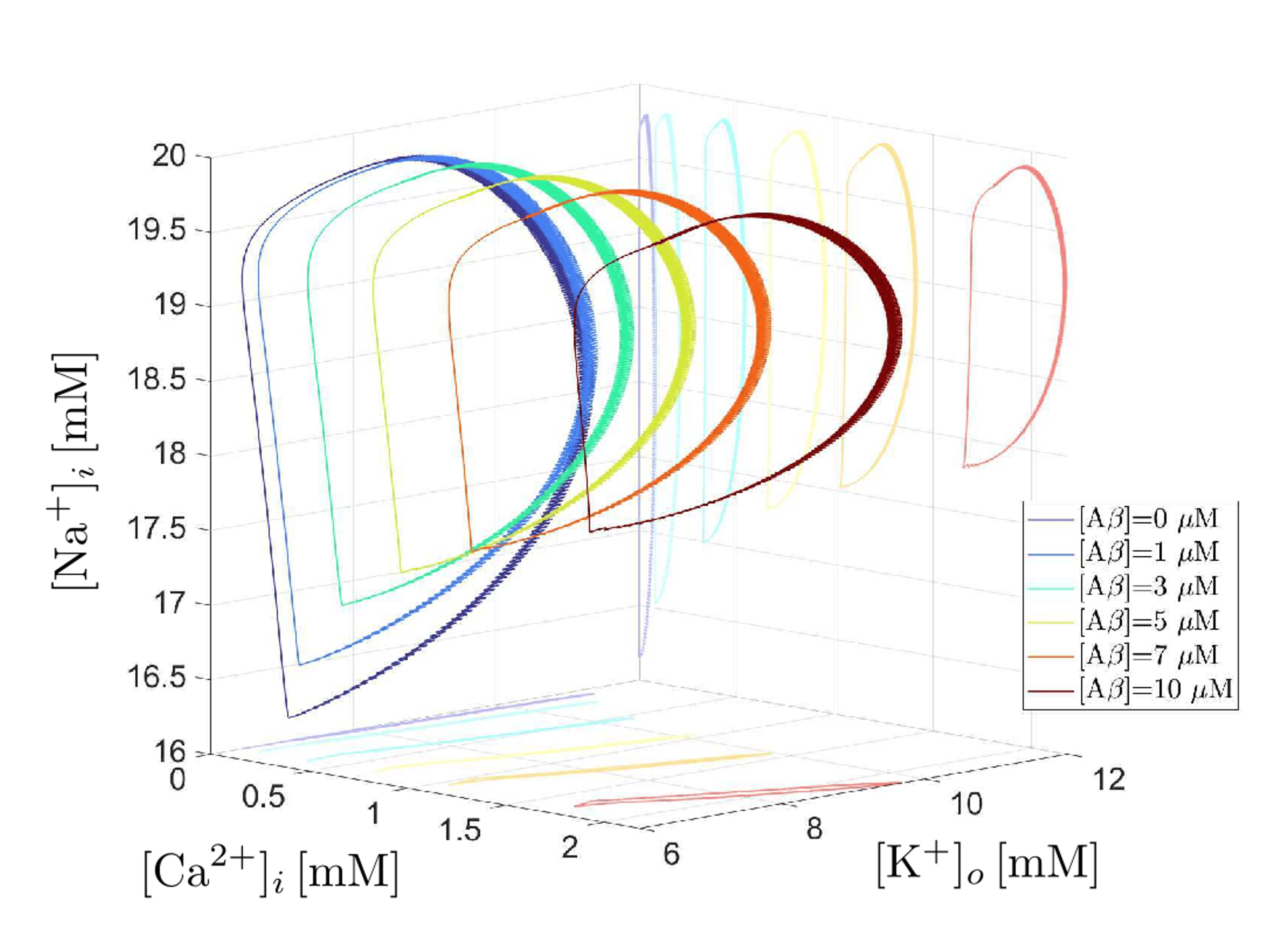}   
     \caption{Stable attractor w.r.t. \Abeta{ }concentration}
     \label{fig:attractor_sens}
    \end{subfigure}%
    \begin{subfigure}[t]{0.50\textwidth}
          \includegraphics[width=\textwidth]{NaKnew.png}
          \caption{Stable attractor for $\Ab=10\,\mu\mathrm{M}$ with three details zoomed.}
        \label{fig:attractor_u}
    \end{subfigure}%
    \\
    \centering
    \begin{subfigure}[b]{\textwidth}
        \centering
        \includegraphics[width=0.65\textwidth]{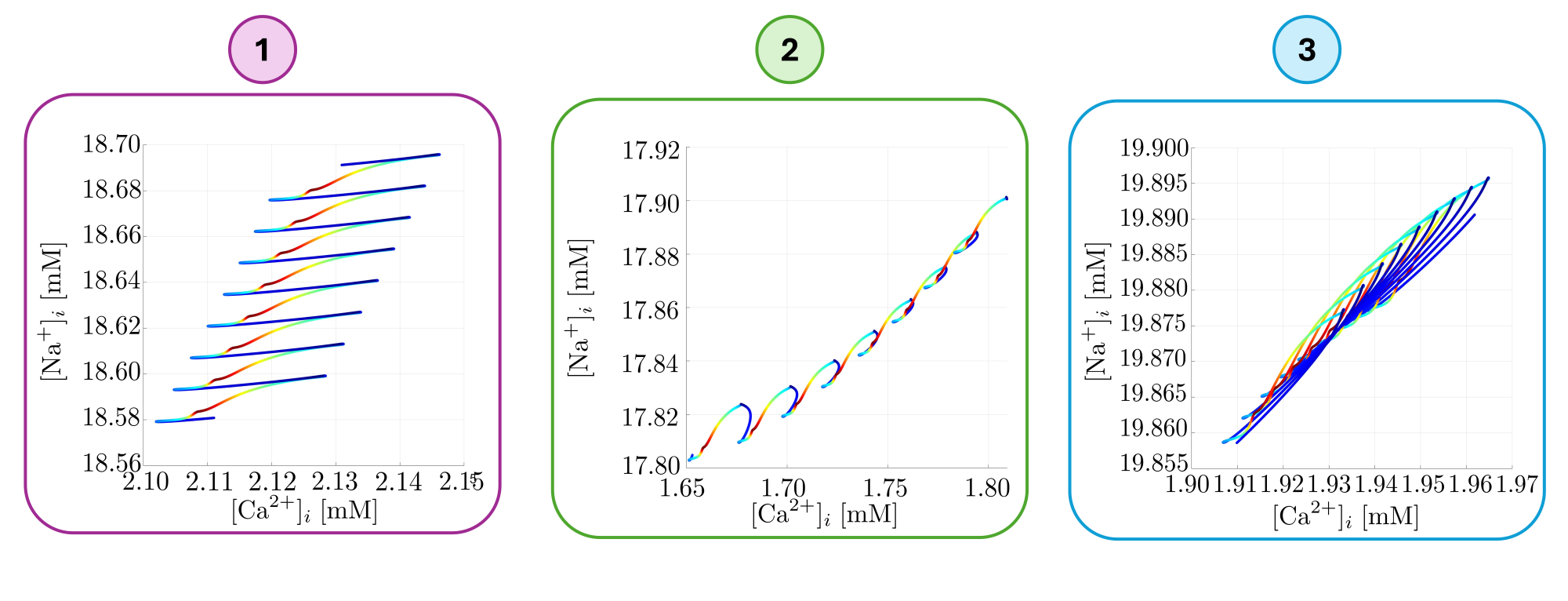}
          \caption{Three zoomed details of stable attractor for $\Ab=10\,\mu\mathrm{M}$.}
        \label{fig:attractor_u_zoom}
    \end{subfigure}%
    \caption{Stable attractor of the ODE system in the 3D space $(\Ca,\K,\Na)$. Sensitivity plot of the attractor with the projection on the three planes (a), attractor colored with $u$-values for $\Ab=10\,\mu\mathrm{M}$ (b), and three details zoomed in the $(\Ca,\Na)$ plane (c).}
    \label{fig:attractor}
\end{figure}
Figure~\ref{fig:attractor} illustrates the stable attractor of the ODE system describing the coupled ionic dynamics involving intracellular calcium $\Ca$, extracellular potassium $K$, and intracellular sodium $\Na$, which are the three critical concentrations for regulating neuronal excitability and transmembrane potential behavior. The attractor shows a toroidal shape associated with a quasi-periodic system dynamics, observed also in the classical Barreto-Cressman ionic model for epileptic seizure modeling \cite{erhardt_dynamics_2020}. Figure \ref{fig:attractor_sens} shows the attractor projected in three-dimensional space for different \Abeta{} concentration values. As \Abeta{} increases, the shape and position of the attractor change noticeably. Specifically, the attractor undergoes a flattening and shifting along the $\Ca$ axis, which reflects a progressive alteration in the homeostatic balance of the system. In particular, we can notice a sensible increase in calcium concentration \cite{berridge_calcium_2010}. These changes suggest that higher levels of \Abeta{} affect the intrinsic dynamics of ion regulation, potentially leading to altered excitability or a transition toward pathological activity patterns. 
\par
Figure \ref{fig:attractor_u} shows the attractor obtained for  $\Ab = 10\,\mu\mathrm{M}$, with the trajectory colored according to the transmembrane potential $u$. In Figure \ref{fig:attractor_u_zoom}, we highlight three different views the attractor, labeled (1), (2), and (3) in the boxes. (1) shows a trajectory that exhibits strong curvature. This indicates a region of dynamical sensitivity, where small changes in ionic variables can produce sensible shifts in the evolution of the system. The box (2) focuses on a segment where the trajectory moves with smoother changes in both potential and ion concentrations. In this region the system evolves into a more stable path. Finally, region (3) shows a portion of the attractor where we have a transition from an active to a quiescent regime in the ionic configuration. The dependence on \Abeta{} concentration supports the idea that the physiological behavior of the neuron is disrupted by the presence of the toxic proteins, with an emergence of an epileptic dynamics \cite{kuchibhotla_abeta_2008}.

\section{Epileptic seizure simulations and interaction with \Abeta{-}aggregates in idealized two-dimensional domains}
\label{sec:2D}
In this section, we present a set of numerical tests with the goal of showing that \Abeta{} can trigger epileptic events at the macro-scale tissue level which propagate throughout the brain tissue. The simulations are based on the monodomain model coupled with the A$\beta$-modified BC model, proposed in Section~\ref{sec:model}. The simulations in this section have been performed using the lymph library \cite{antonietti_lymph_2025}. The numerical simulation of Section 5.3 is run on the GALILEO100 supercomputer (528 computing nodes each 2 x CPU Intel CascadeLake 8260, with 24 cores each, 2.4 GHz, 384GB RAM) at the CINECA supercomputing center.

\subsection{Comparison of different \Abeta{} concentration effects on epileptic seizure propagation}
\label{sec:test_case_square}
The first test case is based on an idealized geometry representing a small portion of grey matter. For this reason, we focus on the isotropic case. The \Abeta{} concentration modifies the neuronal excitability using the model presented in Section \ref{sec:model}. We construct two different test cases to analyze the evolution of the resulting potential values for two values of \Abeta{} ($\Ab = 1\; \mu \mathrm{M} \;$ and $10\; \mu \mathrm{M}$).
\par
\begin{figure}[h]
\centering
\includegraphics[width=\textwidth]{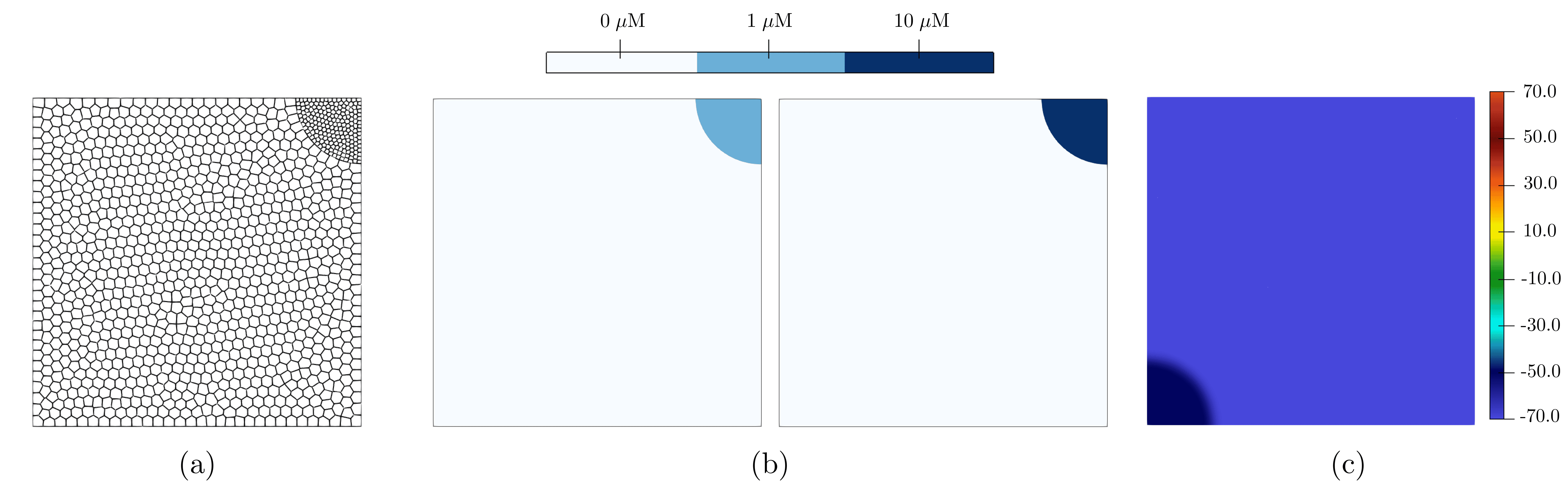}
\caption{(a) Computational mesh; (b) Different pathological values of \Abeta{} concentration; (c) initial condition for the transmembrane potential.} 
\label{fig:square}
\end{figure}

We simulate the evolution of the transmembrane potential in an idealized two-dimensional square domain $\Omega$ of size $(0\,\mathrm{cm},1\,\mathrm{cm})^2$ modelling the grey matter tissue. The domain is divided into two subregions, characterized by different values of \Abeta{}. We localize the pathological \Abeta{} concentration in the subdomain $\OmegaAb= \{ (x,y)\in \Omega\; \mathrm{s.t.}\; (x-1)^2 + (y-1)^2 \leq 0.04 \}$. The complementary region $\Omega \backslash \OmegaAb$ does not contain \Abeta{} (see Figure \ref{fig:square}b). 
We impose an initial localized potential imbalance which models diseased neurons in $\Omega_0$ ($u^0|_{\Omega_0} = -50 \; \mathrm{mV}$). The initial value for the potential is $u^0=-67\;\mathrm{mV}$ in the remaining part of the domain. The initial condition for the transmembrane potential is reported in Figure \ref{fig:square}c. The isotropic conductivity values are taken from \cite{erhardt_dynamics_2020}. 
\par
We report the mesh of $1\,000$ polygonal elements, in Figure \ref{fig:square}a. In particular, we use $800$ elements in $\Omega \backslash \OmegaAb$ and $200$ elements in $\OmegaAb$. The space discretization uses an adaptive PolyDG algorithm with respect to the polynomial degree described in \cite{leimer_saglio_p-adaptive_2025} is used for the discretization. Concerning the time discretization, we consider $\Delta t = 2.5 \, \mu\mathrm{s}$ and $T=200\,\mathrm{ms}$.
\par
\begin{figure}[t!]
    \centering
    \begin{subfigure}[t]{0.423\textwidth}
\includegraphics[width=\textwidth,height=15cm, keepaspectratio]{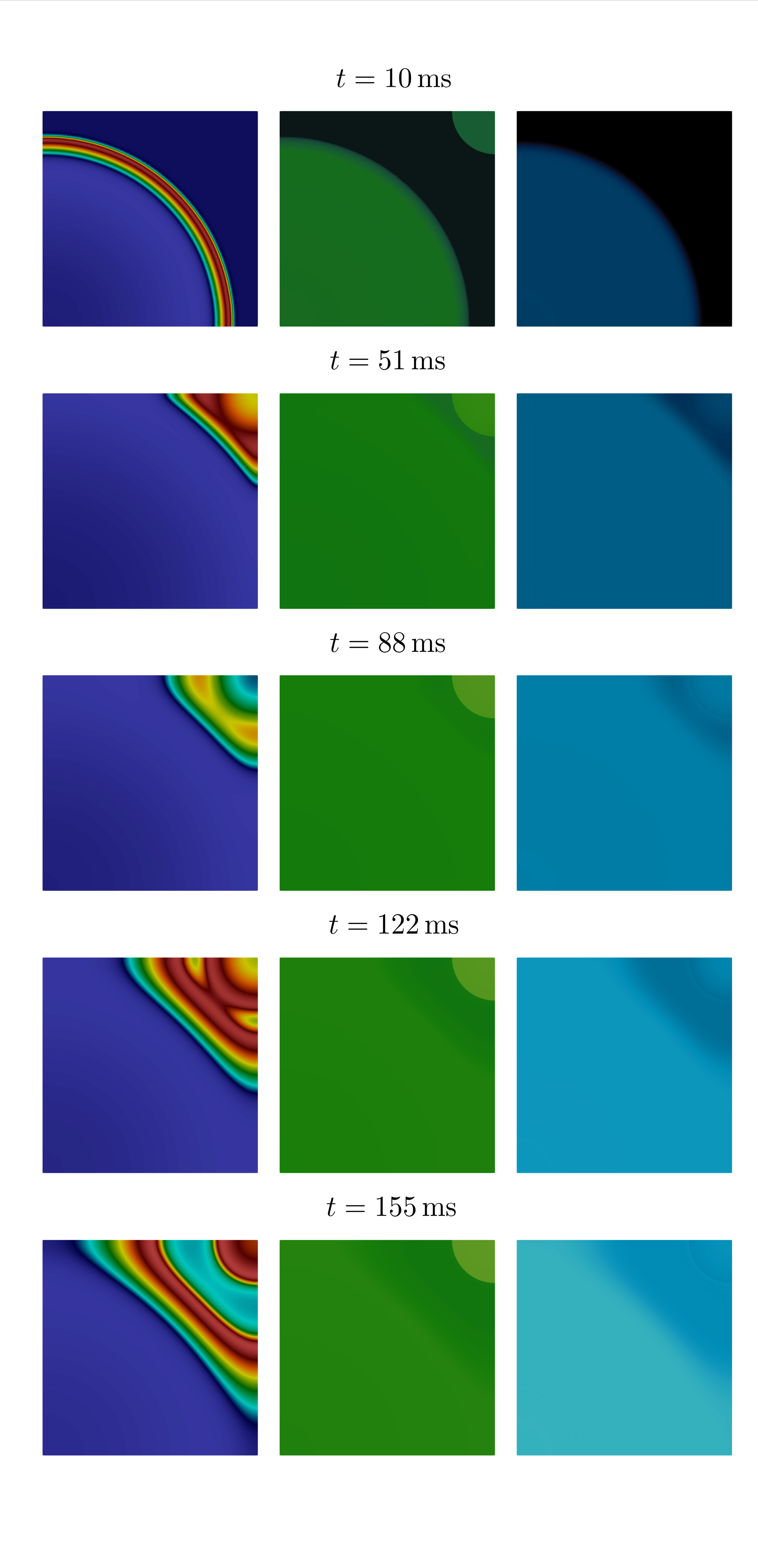}
      \caption{Solution of case $\Ab = 1\, \mu M$: the dynamics mostly governed by the original epileptogenic zone.}
        \label{fig:squaresol:beta1}
    \end{subfigure}\hfill
    \begin{subfigure}[t]{0.555\textwidth}
        \includegraphics[width=\textwidth,height=15cm, keepaspectratio]{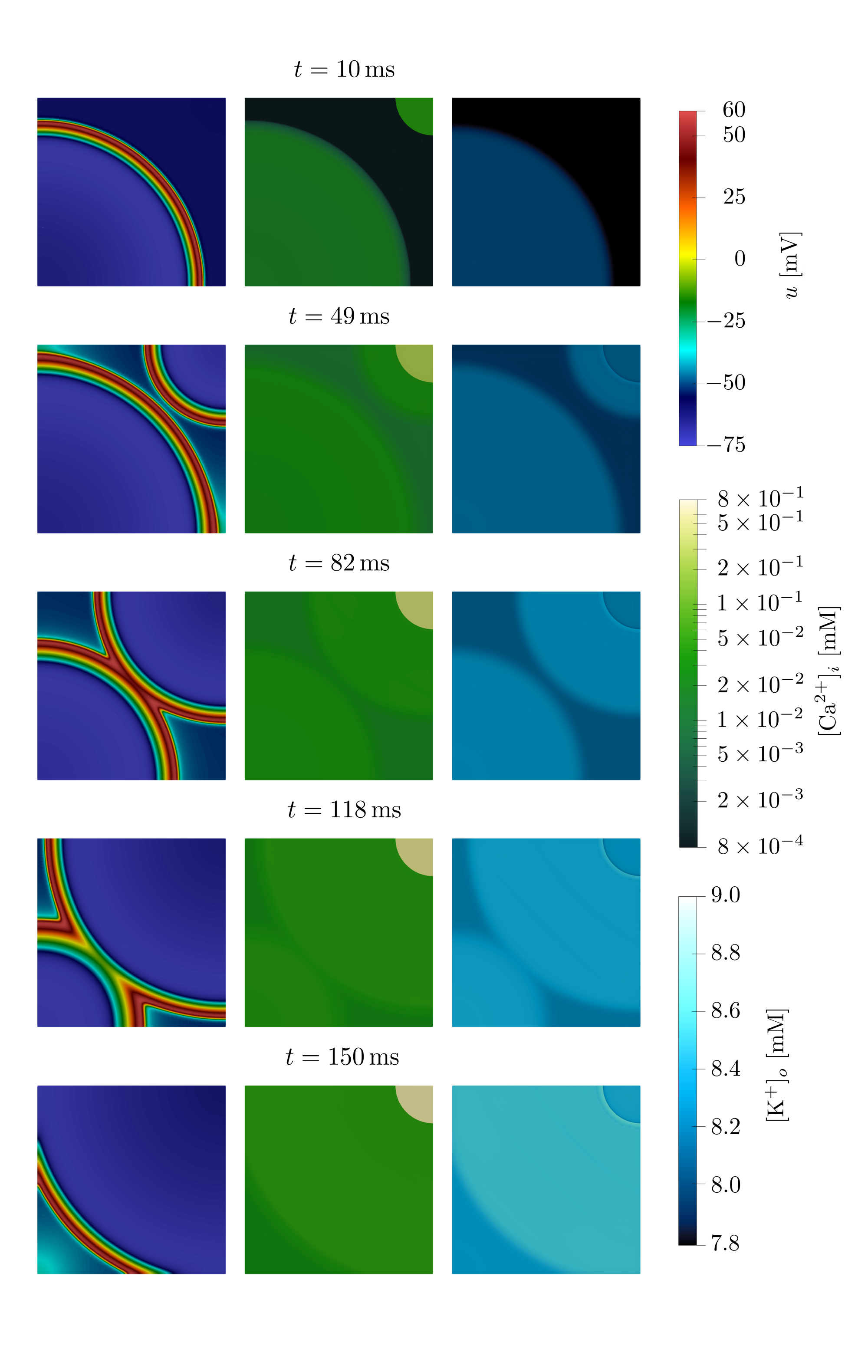}
        \caption{Solution of case $\Ab = 10\, \mu M$: the dominant source of activity shifts to $\OmegaAb$.}
        \label{fig:squaresol:beta10}
    \end{subfigure}
    \caption{Comparison of the evolution of $u$, $\Ca$, $\K$ for $\Ab = 1\,\mu\mathrm{M}$ (a) and $\Ab = 10\,\mu\mathrm{M}$ (b).}
\label{fig:squaresol}
\end{figure}
In Figure \ref{fig:squaresol} we show the evolution of transmembrane potential $u$, intracellular calcium concentration $\Ca$, and extracellular potassium concentration $[K]_0$ at different simulation times in the two cases. Taking into account the medium-low concentration $\Ab = 1 \; \mu \mathrm{M}$ (see Figure~\ref{fig:squaresol:beta10}), we observe an increase of the intracellular calcium concentration in $\OmegaAb$ starting from the first milliseconds of evolution, which tends to increase over time. This pathological behaviour is coherent with the results in Section \ref{sec:0d}. On the contrary, the concentration of extracellular potassium $\K$ increases according to the wavefront profile and then decreases during the repolarisation phase of the tissue region. Analysing the evolution of the transmembrane potential, we note that during the temporal evolution from $\OmegaAb$ pathological wavefronts originate with almost the same periodicity and frequency as those born from $\Omega_0$.
Indeed, the pathological wavefronts arise concurrently with the arrival of the wave that originates in $\Omega_0$, disrupting both waves inside the domain. However, waves arising in $\OmegaAb$ do not develop through the whole domain $\Omega$ because they meet the pathological waves arising from $\Omega_0$.
\par
Considering the case of high concentration of $\Ab  = 10 \mu \mathrm{M}$ (see Figure~\ref{fig:squaresol:beta10}), we observe larger values of intracellular calcium concentration in $\OmegaAb$. In particular, considering the time $t=118 \;\mathrm{ms}$, we find a concentration $\Ca = 4.98\times 10^{-1} \; \mathrm{mM}$ (compared to a value of $\Ca = 8.89\times 10^{-2} \; \mathrm{mM}$ in the first analysis). The evolution of extracellular potassium shows a significant decrease in the peak values during the wave front propagation, accordingly to the 0D model results of Section~\ref{sec:0d}. The behaviour of the transmembrane potential differs significantly from the previous case. Specifically, pathological wavefronts with different period and frequency emerge from the region $\OmegaAb$, rather than from the initial epileptogenic zone $\Omega_0$. In this simulation, the \Abeta{} pathological region becomes dominant in the production of the pathological seizures.
\begin{figure}[t]
    \centering
    \begin{subfigure}[b]{0.25\textwidth}
    \includegraphics[width=\textwidth]{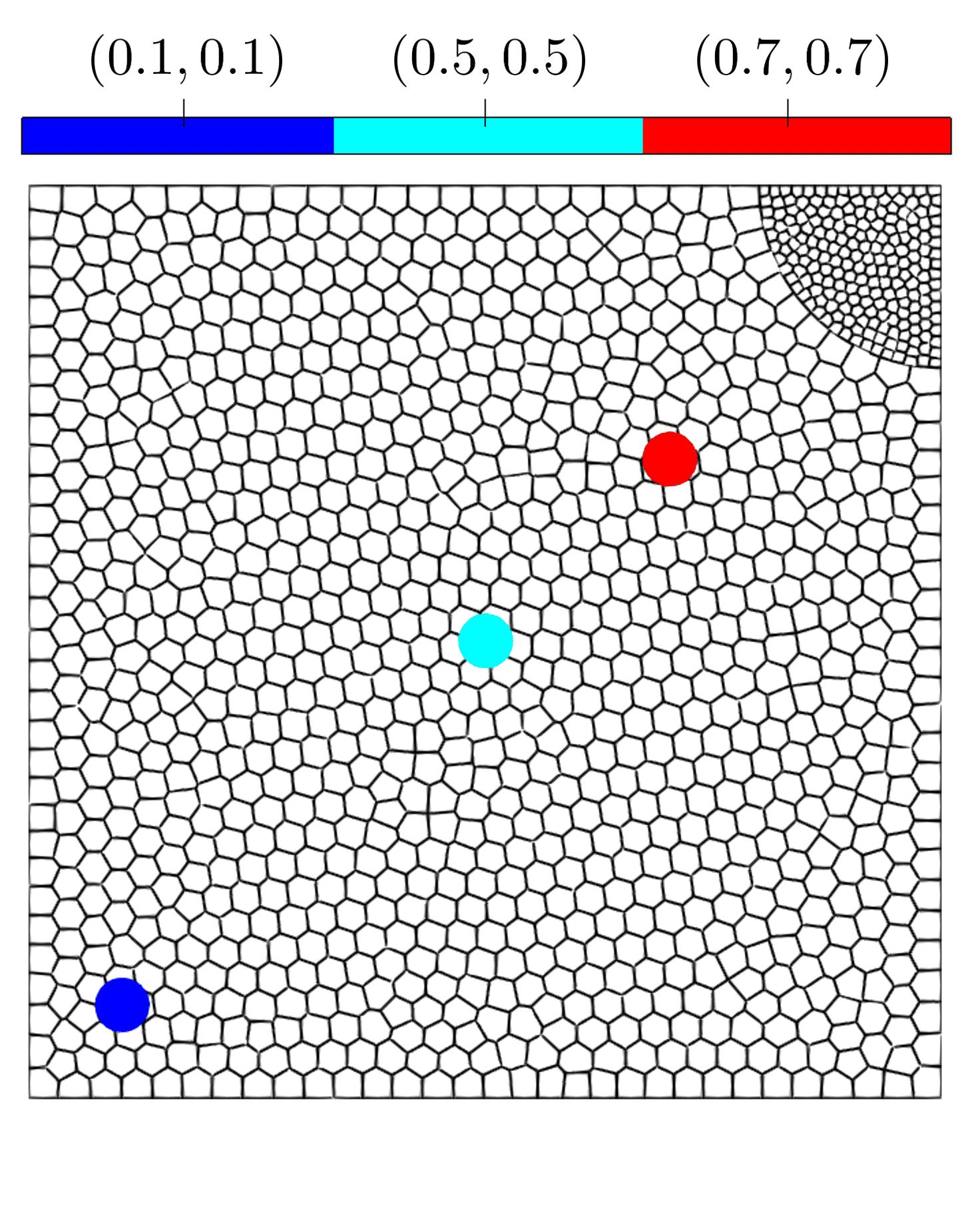}
      \caption{Observation points.}
    \label{fig:potentials:points}
    \end{subfigure}\hfill
    \begin{subfigure}[b]{0.35\textwidth}
    \centering
    \includegraphics[width=\textwidth]{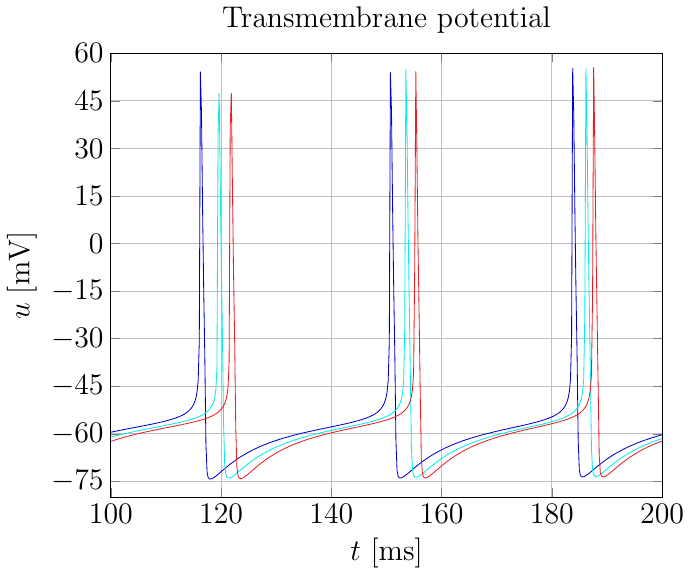}
    \caption{Potential evaluations for $\Ab =  1 \;\mu \mathrm{M}$.}
    \label{fig:potentials:Ab1}
    \end{subfigure}\hfill
    \begin{subfigure}[b]{0.35\textwidth}
    \centering
    \includegraphics[width=\textwidth]{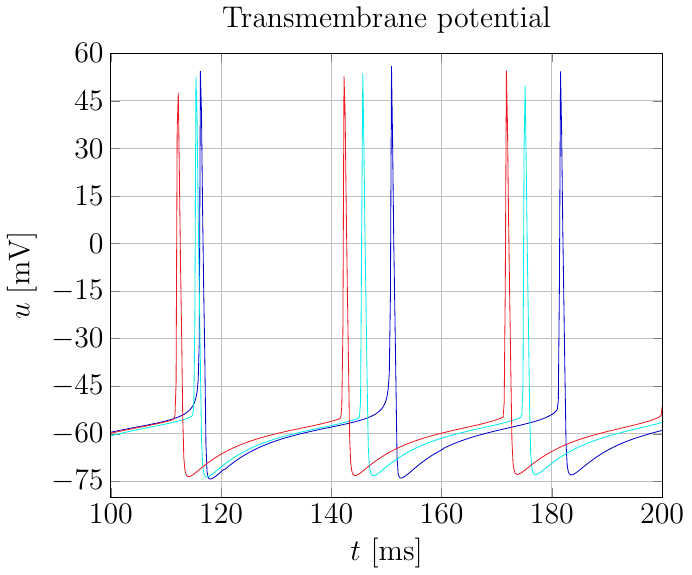}
   \caption{Potential evaluations for $\Ab =  10 \;\mu \mathrm{M}$.}
    \label{fig:potentials:Ab10}
    \end{subfigure}
    \caption{Comparison of transmembrane potential evolution at three spatial points for two concentrations of \Abeta{}. Observation points on the computational mesh (a) are colored to match the curves in the right plots. Moreover, we report the transmembrane potential $u$ over time for $\Ab = 1\;\mu\mathrm{M}$ (b) and $\Ab = 10\;\mu\mathrm{M}$ (c).}
    \label{fig:potentials}
\end{figure}
To illustrate the shift in dynamics, we show the transmembrane potential at three spatial locations for both concentrations of \Abeta{} considered in Figure~\ref{fig:potentials}. The positions of the three points on the computational domain are reported in Figure~\ref{fig:potentials:points}. For low concentrations of \Abeta{}, the system dynamics are governed primarily by the activity originating in the pathological zone in $\Omega_0$. However, when the concentration of \Abeta{} increases to $10\,\mu\mathrm{M}$, the dynamical behaviour changes notably. In particular, Figure~\ref{fig:potentials:Ab10} demonstrates that the red point $(0.7, 0.7)$, located near to $\OmegaAb$, becomes active earlier than the other two analysed points after approximately $100\,\mathrm{ms}$. It is important to note that this change in dynamics also occurs, although more gradually, for lower concentrations of \Abeta{}. As shown in Figure~\ref{fig:potentials:Ab1}, within the first $200\,\mathrm{ms}$, the distance between the activation times reduces sensibly. Nonetheless, the pathological region still acts as an epileptogenic zone, contributing to the overall activity, without overtaking the dynamic control from $\Omega_0$.  This suggests that $\OmegaAb$ acts as a new epileptogenic driver, progressively taking over the system dynamics. 
\par
\begin{figure}[t]
    \centering
    \begin{subfigure}[b]{0.49\textwidth}
    \centering
\includegraphics[width=\textwidth]{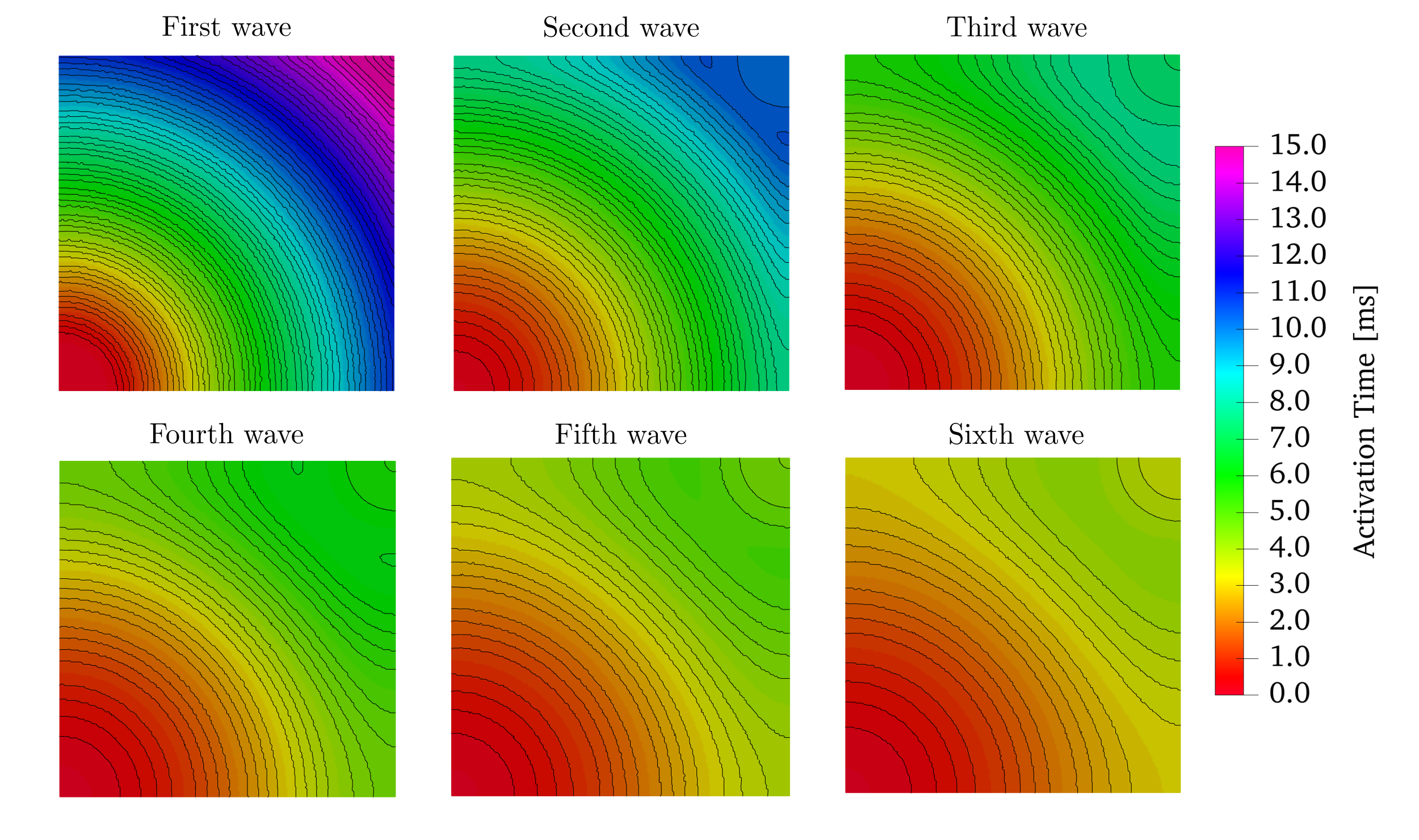}
        \caption{Activation time for $\Ab = 1\,\mu\mathrm{M}$.}
        \label{fig:ATsquare:Ab1}
    \end{subfigure}
    \begin{subfigure}[b]{0.49\textwidth}
    \centering
\includegraphics[width=\textwidth]{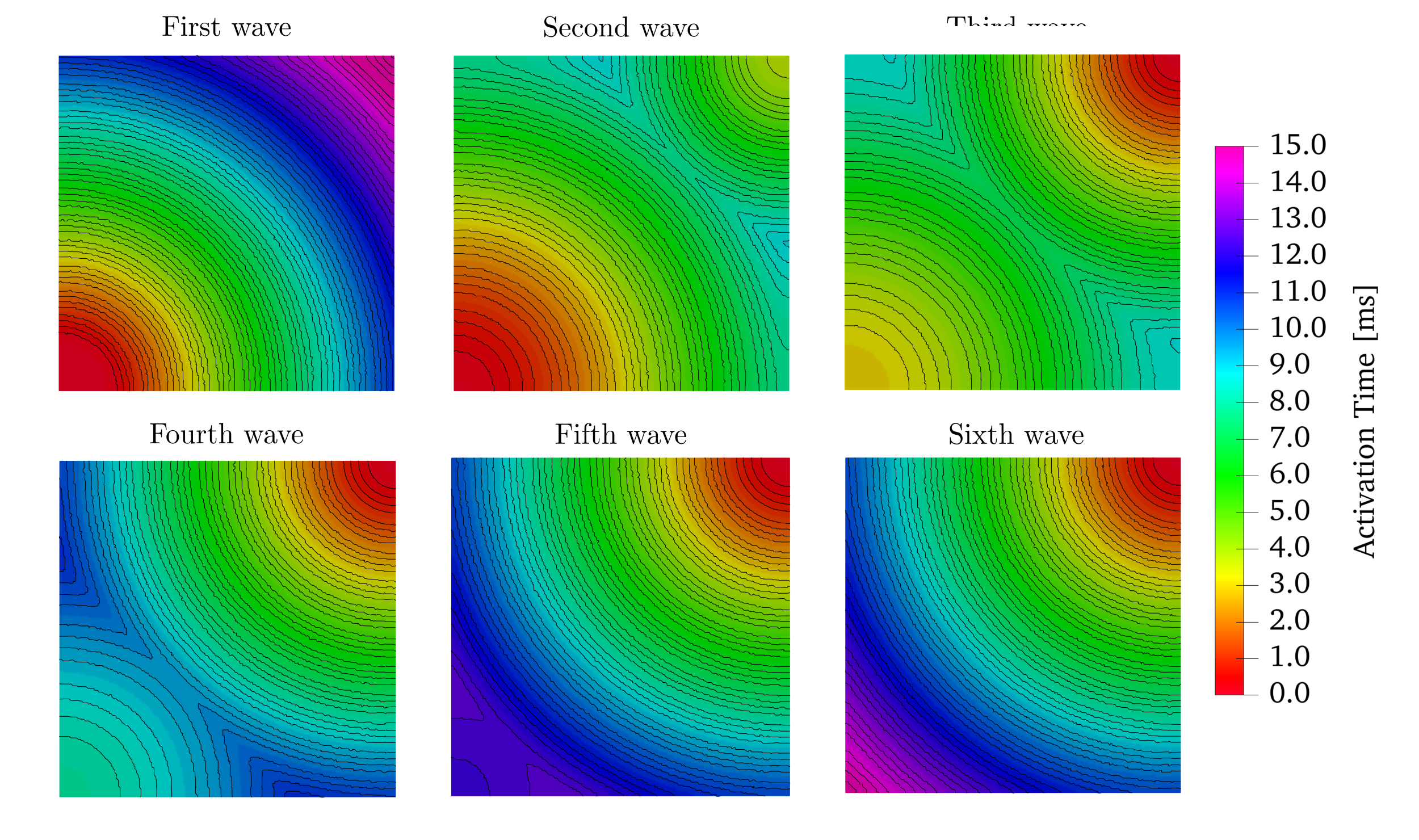}
        \caption{Activation time for $\Ab = 10\,\mu\mathrm{M}$.}
        \label{fig:ATsquare:Ab10}
    \end{subfigure}
    \caption{Maps of activation time of the domain $\Omega$ for the six waves in the two cases: $\Ab = 1\,\mu\mathrm{M}$ (a) and $\Ab = 10\,\mu\mathrm{M}$ (b).} 
    \label{fig:ATsquare}
\end{figure}
Figure~\ref{fig:ATsquare} shows the spatial distribution of activation time \cite{colli_franzone_spreading_1993}, defined for each point $\boldsymbol{x}\in\Omega$ as the first time $\hat{t}$ at which the activation potential is larger than a certain threshold $u_\mathrm{cr}$: 
\begin{equation}
\hat{t} = \underset{t\in[t_\mathrm{min},t_\mathrm{max}]}{\mathrm{argmin}}\left\{u(\boldsymbol{x},t)\geq u_\mathrm{cr}\right\} - t_\mathrm{min}.
\end{equation}
For each wave, we fix $t_\mathrm{min}$ as the first time at which $\exists\, \boldsymbol{x}\in\Omega$ such that $u(\boldsymbol{x},t_\mathrm{min})>0$. In this way, the minimum value of $\hat{t}(\boldsymbol{x})$ is equal to $0$ for all the waves. Each subpanel of Figure~\ref{fig:ATsquare} represents a single activation wave, from the first to the sixth. Figure \ref{fig:ATsquare:Ab1} corresponds to  $\Ab = 1\,\mu\mathrm{M}$. In this case, the activation times show regular and smooth concentric wavefronts. These patterns suggest that the activity is primarily initiated from the original epileptogenic zone and propagates outward in a stable and organized manner. The wavefronts are nearly symmetric, and their shape is preserved across successive activations, indicating minimal disruption in propagation dynamics. Figure \ref{fig:ATsquare:Ab10}, which is associated with $\Ab = 10\,\mu\mathrm{M}$, displays a markedly different behavior. Indeed, the first few waves still exhibit partially concentric fronts; however, the activation patterns become increasingly distorted and asymmetric from the third wave onward. In particular, wavefronts shift spatially and show curvature anomalies, pointing to a secondary or altered activation source. This is consistent with the hypothesis that at higher \Abeta{} concentrations, the pathological region $\OmegaAb$ becomes an active epileptogenic driver \cite{kuchibhotla_abeta_2008}. The time differences also reflect altered excitability and slower or disrupted conduction in certain areas. These observations confirm that the increase in \Abeta{} concentration leads also to a progressive modification of the global dynamics and the wave propagation pattern.
\subsection{Epileptic seizure propagation in presence of multiple \Abeta{} lesions}
\begin{figure}[h]
    \centering
        \includegraphics[width=\textwidth]{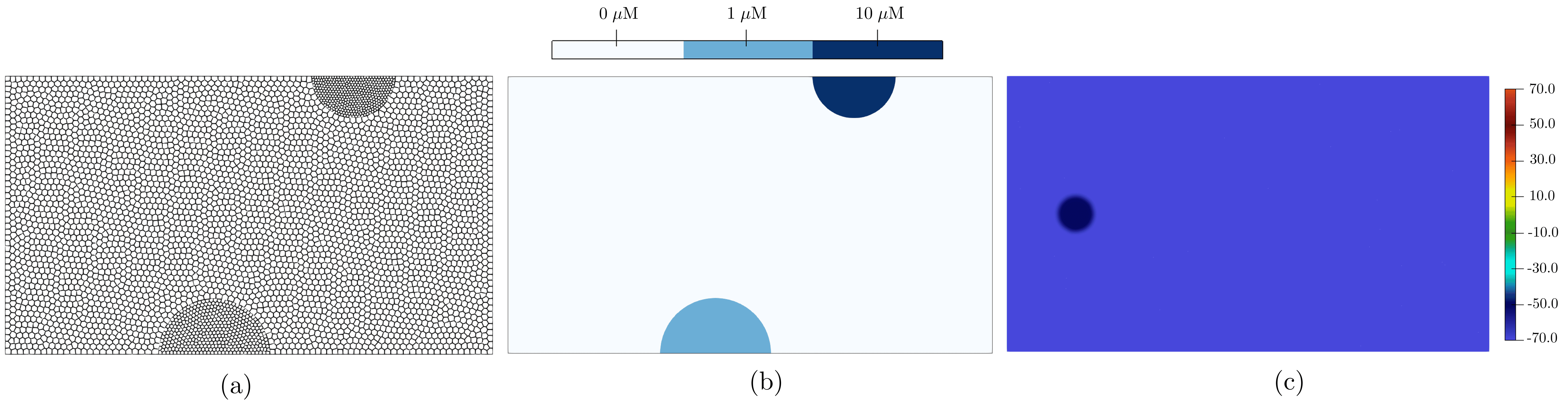}
    \caption{Computational domain: mesh grid (a), with different pathological values of \Abeta{} concentration (b), and the initial condition for the transmembrane potential (c).} 
    \label{fig:secondcomputationaldomain}
\end{figure}

The test case reported in this section is based on an idealized two-dimensional geometry representing a small portion of grey matter tissue, discretized using a polytopal mesh. We consider a domain $\Omega = (-1.5, 2) \times (-1, 1)$, which includes two pathological regions embedded in otherwise healthy tissue.
In particular, two circular region defined as $\Omega_{A\beta}^1 = \left\{ (x, y) \in \Omega \mid (x - 1)^2 + (y - 1)^2 < 0.09\right\}$ and $\Omega_{A\beta}^2 = \left\{ (x, y) \in \Omega \mid x^2 + (y + 1)^2 < 0.16\right\}$ are characterized by an elevated concentration of \Abeta{} peptide, of $10\,\mu\mathrm{M}$ and $1\,\mu\mathrm{M}$, respectively. The complete pathological area $\Omega_{A\beta} = \Omega_{A\beta}^1 \cup \Omega_{A\beta}^2$ is depicted in Figure \hyperref[fig:secondcomputationaldomain]{\ref*{fig:secondcomputationaldomain}b}. A localized subdomain $\Omega_0  = \left\{ (x, y) \in \Omega \mid (x + 1.5)^2 + y^2 < 0.04 \right\}$, located in the left part of the domain, models a cluster of diseased neurons with an initial transmembrane potential set to $u_0|_{\Omega_0} = -50 \, \text{mV}$, thus representing a depolarized and potentially excitable area capable of initiating abnormal signal propagation (see Figure \hyperref[fig:secondcomputationaldomain]{\ref*{fig:secondcomputationaldomain}c}. The remaining part of the domain is associated with an initial transmembrane potential of $ u_0|_{\Omega_h} = -67 \, \text{mV}$.

\begin{figure}[h]
    \centering
    \includegraphics[width=\textwidth]{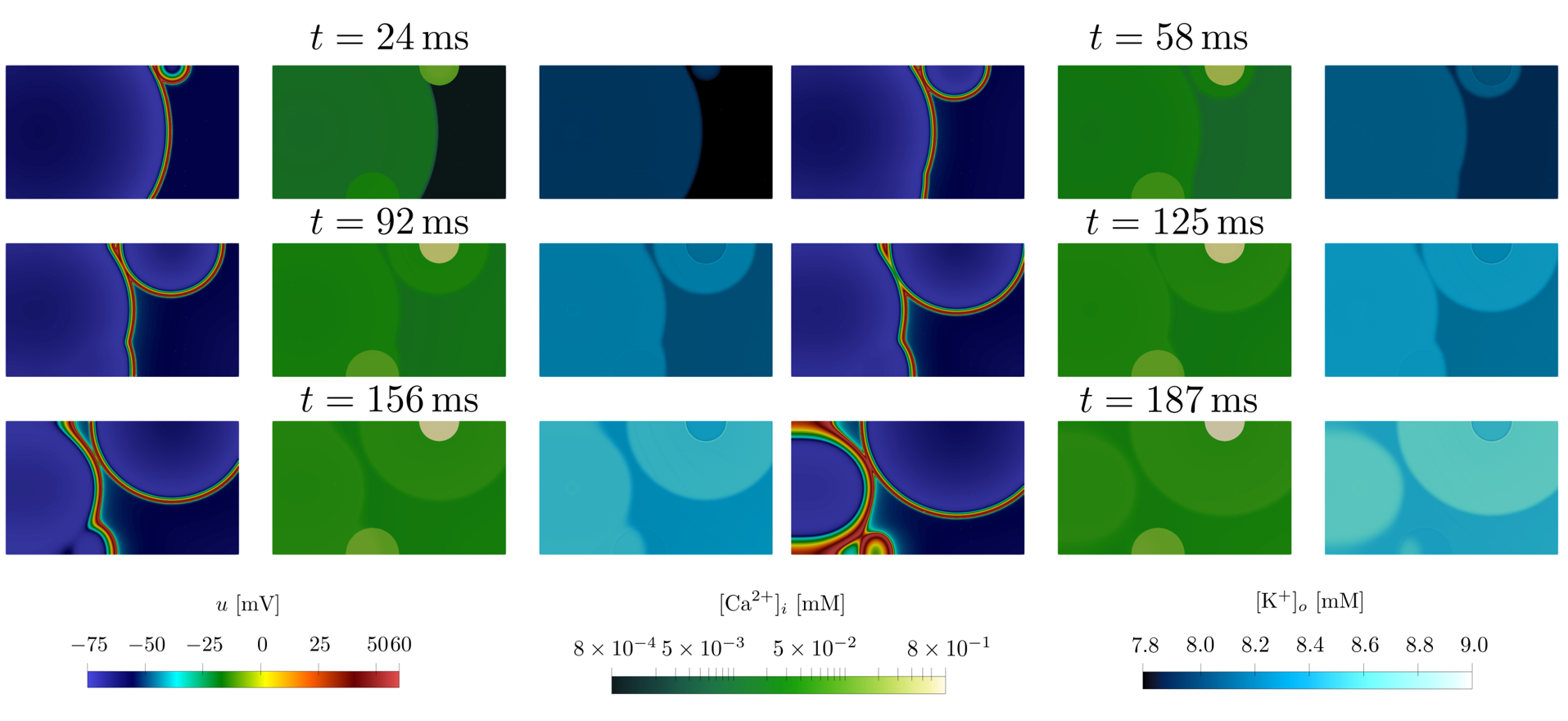}
    \caption{Comparison of the evolution of $u$, $\Ca$, $\K$ for the test case with two different \Abeta{} concentrations inside the computational domain.} 
    \label{fig:BetaDouble}
\end{figure}

Concerning the problem discretization, we report the mesh of $4\,800$ polygonal elements, in Figure \hyperref[fig:secondcomputationaldomain]{\ref*{fig:secondcomputationaldomain}a}. In particular, we use $4000$ elements in $\Omega \backslash \OmegaAb$, $400$ elements in $\OmegaAb^1$ and $400$ elements in $\OmegaAb^2$. The space discretization uses an adaptive PolyDG algorithm with respect to the polynomial degree described in \cite{leimer_saglio_p-adaptive_2025} is used for the discretization. Concerning the time discretization, we consider $\Delta t = 2.5 \, \mu\mathrm{s}$ and $T=200\,\mathrm{ms}$. Figure \ref{fig:BetaDouble} illustrates the simulation results on the computational grid, representing the evolution of the transmembrane potential (left) and both calcium (center) and potassium concentration (right). This scenario builds upon the earlier findings reported in Section \ref{sec:0d}, where increased $\Ab$ levels were shown to enhance intracellular calcium concentration, induce persistent neuronal hyperexcitability, and promote the emergence of pathological spiking activity. 
In line with these observations, the pathological region $\Omega_{\mathrm{A}\beta}$ becomes an independent epileptogenic source. Pathological wavefronts originate from $\Omega_{\mathrm{A}\beta}$ and exhibit different frequency and periodicity compared to the initial epileptogenic zone. This shift in the dominant source of activity underlines the relevance of $\Ab$ on the global propagation dynamics of the system. In particular, it is evident that the stronger dynamics is produced by the larger concentration in the region $\Omega_{A\beta}^1$, unless the size of the domain is smaller than that of $\Omega_{A\beta}^2$. Coherently with the sensitivity analysis performed in the 0D model (Section \ref{sec:0d}), where high $\Ab$ caused an early and sustained elevation of $\Ca$, along with altered $\K$ dynamics. 
\begin{figure}[t]
    \centering
\includegraphics[width=\textwidth]{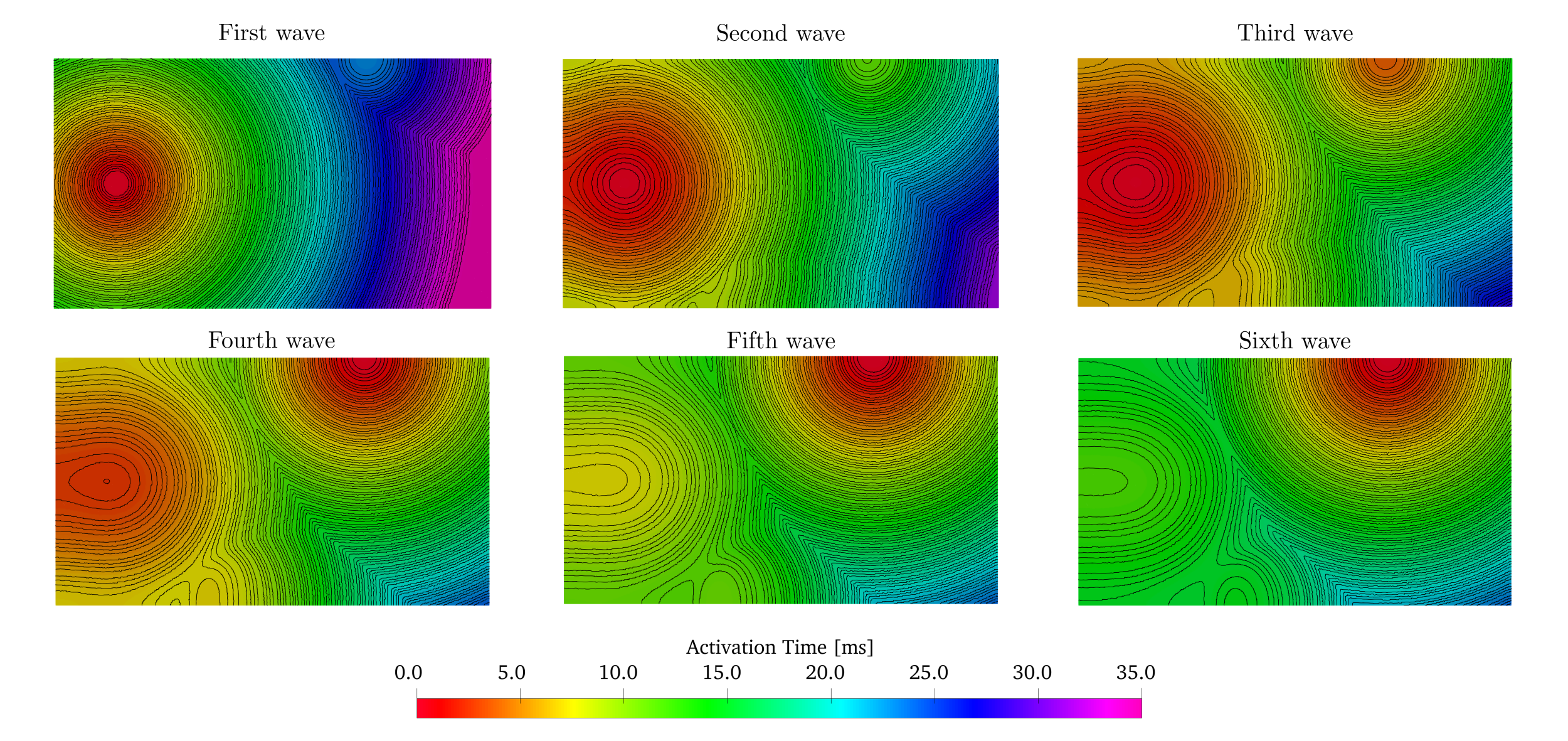}
    \caption{Maps of activation time of the domain $\Omega$ for the six waves.} 
    \label{fig:ATabeta}
\end{figure}
Moreover, the spatiotemporal distributions of extracellular potassium and intracellular calcium further underscore the \Abeta{-}driven pathology: in the vicinity of $\Omega_{\mathrm{A}\beta}$, $\K$ rises sharply at the leading edge of each wavefront, creating a depolarized region that lowers the firing threshold of neighboring neurons, while $\Ca$ exhibits a prolonged plateau with superimposed oscillations that follow after the wave passes, indicative of impaired calcium clearance and repetitive bursting. This pathological setting can dynamically transform a previously healthy region into a dominant epileptogenic driver, thus providing a mechanistic explanation for how regions affected by Alzheimer’s disease can become epileptogenic sources \cite{kuchibhotla_abeta_2008}. Figure \ref{fig:ATabeta} presents the activation times for different waves inside the domain \cite{colli_franzone_spreading_1993}. This result illustrates how the presence of high \Abeta{} levels alters the spatiotemporal propagation of epileptic activity. The simulations show that the region $\Omega_{A\beta}^1$ with high \Abeta{-}concentrations ($10\,\mu\mathrm{M}$) acts as an independent source of excitation, governing the spatiotemporal dynamics against both the initial epileptic region and the second pathological amyloid region $\Omega_{A\beta}^2$. This creates asynchronous and heterogeneous wavefronts, causing a disruption of the regularity of the propagation pattern.
\subsection{Epileptic seizure propagation on realistic brain sections with PET-derived \Abeta{-}distributions}
\label{sec:realistic}

In this section, we present a numerical simulation performed on a realistic two-dimensional brain slice geometry, where the distribution of A\textbeta{} peptides is derived from positron emission tomography (PET) imaging data. 
\begin{figure}[h]
    \centering
    \includegraphics[width=\textwidth]{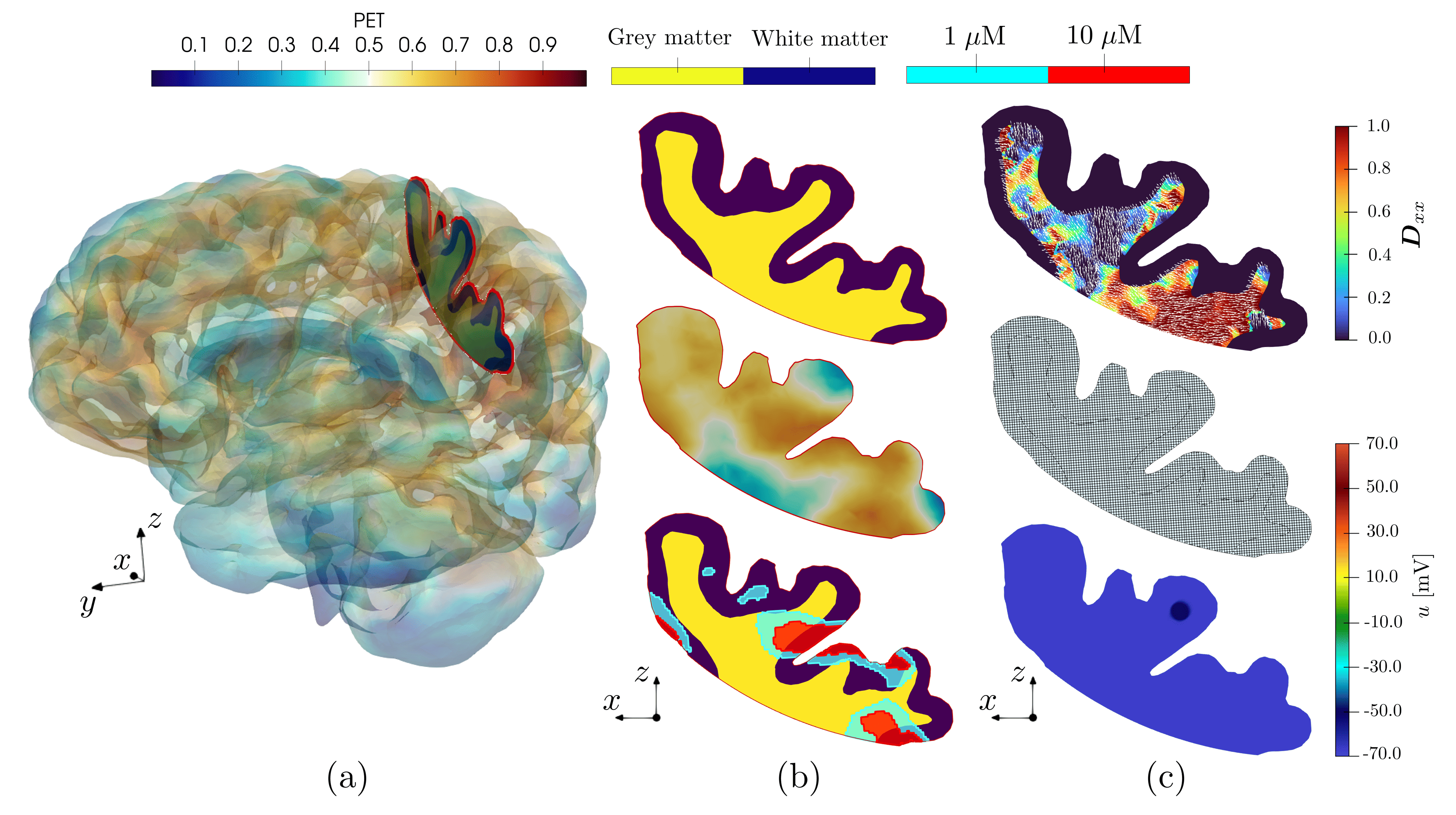}
    \caption{Setup of the realistic brain simulation: (a) Three-dimensional distribution of PET concentration in the human brain. (b) Top: Sagittal slices showing the distinction between grey matter (yellow) and white matter (blue); center: PET distribution in selected section; bottom: \Abeta{} distribution in the computational domain. (c) Top: $\boldsymbol{D}_{xx}$ anisotropic axonal directions for white matter tissue; center: computational domain with polytopal mesh of 6532 elements ($h = 0.10201$); bottom: initial pathological condition for the transmembrane potential.}
    \label{fig:Brainsetup}
\end{figure}
We consider the subject OAS30080 (see Figure \ref{fig:Brainsetup}a), in the OASIS dataset \cite{lamontagne_oasis_2019}, for which the geometry is a portion of a section on the coronal plane segmented from a structural magnetic resonance image (MRI) by means of Freesurfer \cite{fischl_freesurfer_2012} (see Figure \ref{fig:Brainsetup}b). We derive the patterns of A\textbeta{} concentration from PET images with Pittsburgh Compound-B to feed our ionic model in equation \eqref{eq:bc_abeta}. We analyze a pathological configuration, corresponding to a late-stage evolution of Alzheimer’s disease, with elevated and widespread A\textbeta{} concentrations, extending into subcortical regions. We consider a concentration of $[A\beta] = 1\,\mu\mathrm{M}$ for PET values inside the interval $(0.65,0.70)$ and $[A\beta] = 10\,\mu\mathrm{M}$ for PET values higher than $0.70$, as shown in Figure \ref{fig:Brainsetup}b.
In the white matter, we assume that the conduction tensor $\mathbf{\Sigma}$ has the following structure:
\begin{equation}
\label{eq:sigma_tensor}
    \mathbf{\Sigma}(\boldsymbol{x}) = \sigma_\mathrm{iso}\mathbf{I} + \sigma_\mathrm{axn}(\boldsymbol{a}(\boldsymbol{x})\otimes \boldsymbol{a}(\boldsymbol{x})),
\end{equation}
where $\boldsymbol{a}=\boldsymbol{a}(\boldsymbol{x})$ is the direction of axonal fibres in $\boldsymbol{x}\in\Omega$ and the conduction is associated with an isotropic component $\sigma_\mathrm{iso}=0.0735\,\mathrm{S\,m}^{-1}$ and an axonal one $ \sigma_\mathrm{axn}=0.6\,\mathrm{S\,m}^{-1}$. The axonal direction is derived from diffusion weighted imaging (DWI) using Freesurfer \cite{fischl_freesurfer_2012} (see Figure \ref{fig:Brainsetup}c).
\par
The domain is discretized using a polytopal mesh constructed starting from a Cartesian mesh, where the elements are then cut in correspondence with the boundary and the internal grey-white matter interfaces. This approach allows us to retain the geometric complexity of the brain structures while preserving the efficiency of structured discretizations, exploiting the PolyDG formulation described in Section \ref{sec:3}. The final grid is composed of $6\,532$ polygonal elements (see Figure \ref{fig:Brainsetup}c). An initial localized perturbation of the transmembrane potential is imposed to trigger epileptic activity, simulating an epileptogenic focus, represented in Figure \ref{fig:Brainsetup}c. 
 \begin{figure}[h]
    \centering
    \includegraphics[width=\textwidth]{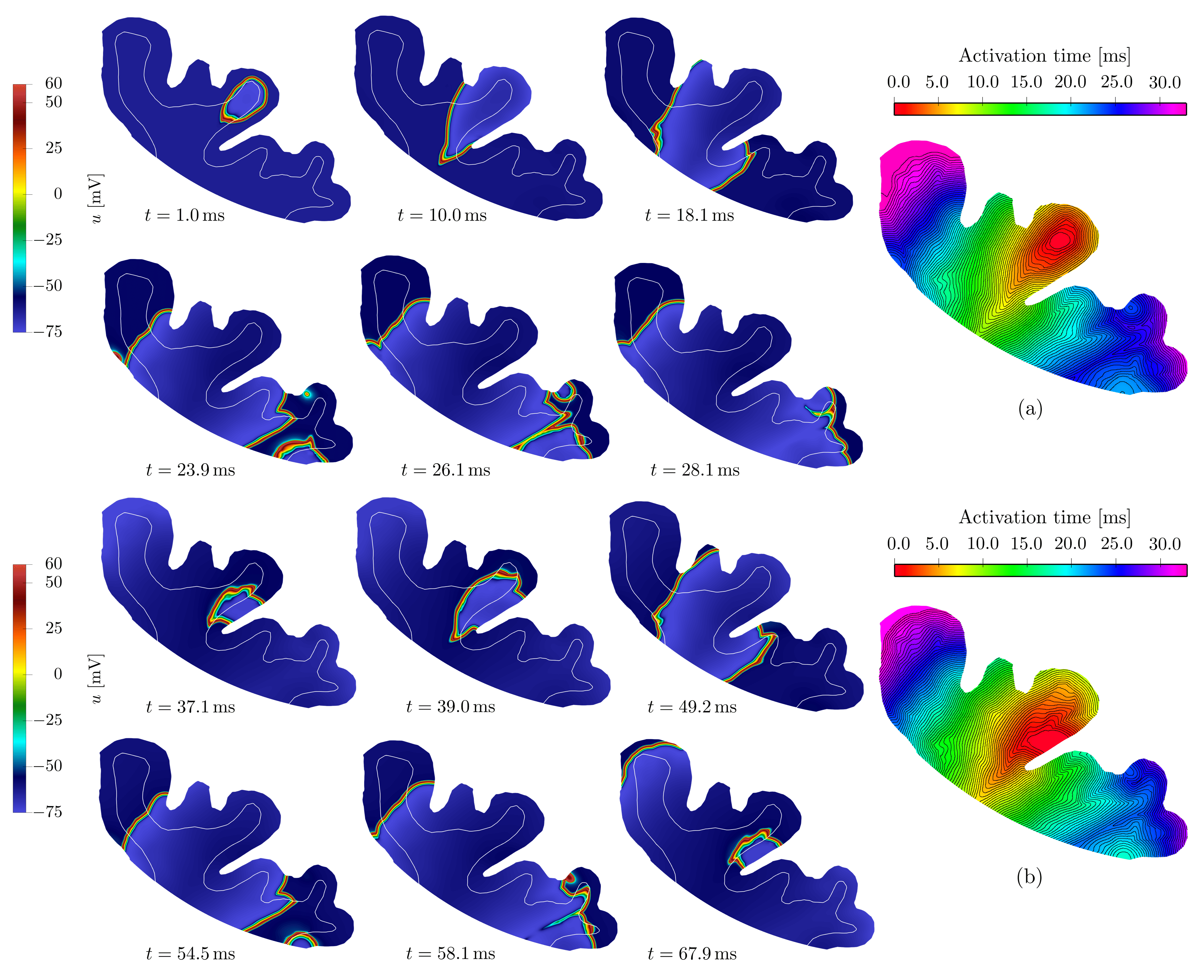}
    \caption{Temporal evolution of the transmembrane potential $u$ for different time instants and corresponding activation time maps for the two waves: (a) first wave (b) second wave. The color scale represents the membrane potential $u$ (in mV), while contour lines highlight the propagating wavefronts across the tissue domain.} 
    \label{fig:BetaBrain}
\end{figure}
\par
Figure~\ref{fig:BetaBrain} illustrates the evolution of the transmembrane potential $u$ at different time instants and the corresponding activation time maps. In this configuration, the wavefront starts at the epileptogenic zone and propagates through the tissue. 
In the isotropic grey matter region, the transmembrane potential evolution remains concentric as expected because of the absence of structural preferential directions. In contrast, the anisotropic conductivity tensor creates preferential propagation along the main axonal directions in the white matter, stretching the wavefronts along these preferential directions.
\par
The heterogeneous distribution of A\textbeta{} further modulates the spatiotemporal dynamics of the electrical activity: regions characterized by higher values of A\textbeta{} exhibit altered excitability and act as secondary drivers, locally modifying both the amplitude and timing of the depolarization waves. As a result, the global propagation pattern becomes spatially heterogeneous, reflecting the combined effect of tissue anisotropy and pathological inhomogeneities in  A\textbeta{} concentration.
The depolarization waves exhibit irregular shapes, and spatially shifted activation centers are present, indicating a secondary source of excitation induced by A\textbeta{} accumulation.
These results confirm that increasing A\textbeta{} alters the spatiotemporal organization of the seizure propagation, leading to heterogeneous activation patterns and slower conduction in the affected areas, in agreement with the sensitivity analysis discussed in Section~\ref{sec:0d}.
\par
Conversely, the regions characterized by a moderate \Abeta{} concentration ($\simeq 1\,\mu\mathrm{M}$) are not able to sustain an autonomous firing. The excitatory waves produced by those lesions are rapidly absorbed by the seizures originated by high \Abeta{} regions ($\simeq 10\,\mu\mathrm{M}$), which display lower intrinsic frequency but a higher persistence. The maps of activation time confirm this interplay, revealing a shift of activation centers in space and an irregular depolarization timing associated with high \Abeta{} areas (see Figure~\ref{fig:BetaBrain}). These findings indicate that the seizures are affected by the presence of \Abeta{} lesions both in terms of spatial organization and temporal coherence, making them additional drivers of epilepsy.

\section{Conclusions}
\label{sec:conclusion}
In this work, we have presented a novel computational framework that integrates the effects of amyloid-\textbeta{} accumulation into a detailed ionic model of neuronal electrical activity, extending the Barreto--Cressman formulation for epileptic bursting. We introduced new amyloid-\textbeta{} dependent pathways affecting calcium homeostasis and potassium conductances. These modifications include the inhibition of the plasma membrane $\Capl$-ATPase, the formation of $\Capl$-permeable pores, the overactivation of L-type calcium channels, and the reduction of $\Capl$-sensitive potassium currents.  
The resulting ionic model captures the progressive transition from physiological to pathological excitability as the amyloid-\textbeta{} concentration increases, reproducing key features observed in Alzheimer's disease such as calcium overload, prolonged depolarization, and enhanced burst frequency. The sensitivity analysis confirms that increasing amyloid-\textbeta{} levels drive a transition from quasi-periodic to hyperexcitable dynamics, with sustained elevations of intracellular calcium and reduced potassium peaks, both hallmark signatures of epileptiform behavior. Coupling the modified ionic model with the monodomain formulation enabled the investigation of spatio-temporal seizure propagation in brain tissue. Through the discontinuous Galerkin discretization on polygonal meshes, we simulated both idealized and realistic geometries, including PET-derived amyloid-\textbeta{} distributions. Overall, the proposed model provides a connection between amyloid-\textbeta{} pathology, calcium dysregulation, and epileptic seizure generation in Alzheimer's disease. In particular, the simulations show that tissue regions with high \Abeta{} deposition can themselves act as additional epileptogenic drivers, offering a mechanistic explanation for the increased incidence of seizure-like events observed in patients with Alzheimer’s disease. Future developments will include model calibration with patient-specific imaging data and the extension to three-dimensional brain geometries.

\section*{Acknowledgments}
OASIS-3 provided the brain MRI images: Longitudinal Multimodal Neuroimaging: Principal Investigators: T. Benzinger, D. Marcus, J. Morris; NIH P30 AG066444, P50 AG00561, P30 NS09857781, P01 AG026276, P01 AG003991, R01 AG043434, UL1 TR000448, R01 EB009352. AV-45 doses were provided by Avid Radiopharmaceuticals, a wholly-owned subsidiary of Eli Lilly.

\bibliographystyle{hieeetr}
\bibliography{sample.bib}

@article{schreiner_simulating_2022,
  title     = {Simulating epileptic seizures using the bidomain model},
  author    = {Schreiner, Jakob and Mardal, Kent-Andre},
  journal   = {Scientific Reports},
  year      = {2022},
  volume    = {12},
  number    = {1},
  pages     = {10065},
  doi       = {10.1038/s41598-022-12101-y},
}

@article{erhardt_dynamics_2020,
    title       = {Dynamics of a neuron--glia system: the occurrence of seizures and the influence of electroconvulsive stimuli: A mathematical and numerical study},
    author      = {Erhardt, Andr{\'e} H and Mardal, Kent-Andre and Schreiner, Jakob E},
    journal     = {Journal of Computational Neuroscience},
    year        = {2020},
    volume      = {48},
    pages       = {229--251},
    doi         = {10.1007/s10827-020-00746-5},
}

@article{cressman_influence_2009,
    title       = {The influence of sodium and potassium dynamics on excitability, seizure, and the stability of persistent states: {I}. {S}ingle neuron dynamics},
    author      = {Cressman, J. R. and Ullah, G. and Ziburkus, J. and Schiff, S. J. and Barreto, E.},
    journal     = {Journal of Computational Neuroscience},
    year        = {2009},
    volume      = {26},
    pages       = {159--170},
    doi         = {10.1007/s10827-008-0132-4}
}

@article{ullah_influence_2009,
    title       = {The influence of sodium and potassium dynamics on excitability, seizures, and the stability of persistent states: II. {N}etwork and glial dynamics},
    author      = {Ullah, Ghanim and Cressman Jr., John R and Barreto, Ernest and Schiff, Steven J},
    journal     = {Journal of Computational Neuroscience},
    year        = {2009},
    volume      = {26},
    pages       = {171--183},
    doi         = {10.1007/s10827-008-0132-4}
}

@article{barreto_ion_2011,
    title       = {Ion concentration dynamics as a mechanism for neuronal bursting},
    author      = {Barreto, Ernest and Cressman, John R},
    journal     = {Journal of Biological Physics},
    year        = {2011},
    volume      = {37},
    pages       = {361--373},
    doi         = {10.1007/s10867-010-9212-6},
}

@article{steinlein_calcium_2014,
    title       = {Calcium signaling and epilepsy},
    author      = {Steinlein, O. K.},
    journal     = {Cell and Tissue Research},
    year        = {2014},
    volume      = {357},
    pages       = {385--393},
    doi         = {10.1007/s00441-014-1849-1}
}

@article{kuchibhotla_abeta_2008,
    title        = {A$\beta$ plaques lead to aberrant regulation of calcium homeostasis in vivo resulting in structural and functional disruption of neuronal networks},
    author       = {Kuchibhotla, K. V. and Goldman, S. T. and Lattarulo, C. R. and Wu, H. and Hyman, B. T. and Bacskai, B. J.},
    journal      = {Neuron},
    year         = {2008},
    volume       = {59},
    number       = {2},
    pages        = {214--225},
    doi          = {10.1016/j.neuron.2008.06.008},
}

@article{berridge_calcium_2010,
    title        = {Calcium hypothesis of  {A}lzheimer’s disease},
    author       = {Berridge, M. J.},
    journal      = {Pflügers Archiv - European Journal of Physiology},
    year         = {2010},
    volume       = {459},
    pages        = {441--449},
    doi          = {10.1007/s00424-009-0736-1}
}

@article{vonbonhorst_impact_2022,
    title       = {Impact of $\beta$-amyloids induced disruption of $\mathrm{Ca}^{2+}$ homeostasis in a simple model of neuronal activity},
    author      = {Prista von Bonhorst, F. and Gall, D. and Dupont, G.},
    journal     = {Cells},
    year        = {2022},
    volume      = {11},
    pages       = {615},
    doi         = {10.3390/cells11040615},
}

@article{latulippe_mathematical_2018,
    title       = {A mathematical model for the effects of amyloid beta on intracellular calcium},
    author      = {Latulippe, J. and Lotito, D. and Murby, D.},
    journal     = {PLOS One},
    year        = {2018},
    volume      = {13},
    number      = {8},
    pages       = {1-27},
    doi         = {10.1371/journal.pone.0202503},
}

@article{raskatov_what_2019,
    title       = {What is the "relevant" amyloid $\beta$42 concentration?},
    author      = {Raskatov, J. A.},
    journal     = {Chembiochem},
    year        = {2019},
    volume      = {20},
    number      = {13},
    pages       = {1725--1726},
    doi         = {10.1002/cbic.201900097},
}

@article{berrocal_calmodulin_2012,
    title       = {Calmodulin antagonizes amyloid-$\beta$ peptides-mediated inhibition of brain plasma membrane $\mathrm{{C}a}^{2+}$-{ATP}ase},
    author      = {Berrocal, M. and Sep{\'u}lveda, M. R. and V{\'a}zquez-Hern{\'a}ndez, M. and Mata, A. M.},
    journal     = {Biochimica et Biophysica Acta (BBA) - Molecular Basis of Disease},
    year        = {2012},
    volume      = {1822},
    pages       = {961--969},
    doi         = {10.1016/j.bbadis.2012.02.013},
}

@article{berrocal_altered_2009,
    title       = {Altered $\mathrm{{C}a}^{2+}$ dependence of synaptosomal plasma membrane $\mathrm{{C}a}^{2+}$-{ATP}ase in human brain affected by {A}lzheimer's disease},
    author      = {Berrocal, M. and Marcos, D. and Sepúlveda, M. R. and Pérez, M. and Ávila, J. and Mata, A. M.},
    journal     = {The FASEB Journal},
    year        = {2009},
    volume      = {23},
    pages       = {1826--1834},
    doi         = {10.1096/fj.08-121459}
}

@article{mirdha_aggregation_2024,
    title       = {Aggregation behavior of amyloid beta peptide depends upon the membrane lipid composition},
    author      = {Mirdha, L.},
    journal     = {The Journal of Membrane Biology},
    year        = {2024},
    volume      = {257},
    pages       = {151--164},
    doi         = {10.1007/s00232-024-00314-3},
}

@article{good_effect_1996,
    title       = {Effect of $\beta$-amyloid block of the fast-inactivating $\mathrm{K}^{+}$ channel on intracellular $\mathrm{Ca}^{2+}$ and excitability in a modeled neuron},
    author      = {Good, T. A. and Murphy, R. M.},
    journal     = {Proceedings of the National Academy of Sciences},
    year        = {1996},
    volume      = {93},
    number      = {26},
    pages       = {15130--15135},
    doi         = {10.1073/pnas.93.26.15130}
}

@article{ishii_amyloid_2019,
    title       = {Amyloid-beta modulates low-threshold activated voltage-gated {L}-type calcium channels of arcuate neuropeptide {Y} neurons leading to calcium dysregulation and hypothalamic dysfunction},
    author      = {Ishii, M. and Hiller, A. J. and Pham, L. and McGuire, M. J. and Iadecola, C. and Wang, G.},
    journal     = {Journal of Neuroscience},
    year        = {2019},
    volume      = {39},
    number      = {44},
    pages       = {8816--8825},
    doi         = {10.1523/JNEUROSCI.0617-19.2019}
}

@article{kim_effects_2011,
    title       ={Effects of Amyloid-$\beta$ peptides on voltage-gated {L}-type {C}a{V}1.2 and {C}a{V}1.3 $\mathrm{Ca}^{2+}$ channels},
    author      = {Kim, S. and Rhim, H.},
    journal     = {Molecules and Cells},
    year        = {2011},
    volume      = {32},
    pages       = {289--294},
    doi         = {10.1007/s10059-011-0075-x}
}

@article{yamamoto_suppression_2011,
    title       = {Suppression of a neocortical potassium channel activity by intracellular amyloid-$\beta$ and its rescue with homer1a},
    author      = {Yamamoto, K. and Ueta, Y. and Wang, L. and others},
    journal     = {Journal of Neuroscience},
    year        = {2011},
    volume      = {31},
    number      = {31},
    pages       = {11100--11109},
    doi         = {10.1523/JNEUROSCI.6752-10.2011}
}

@article{yamamoto_amyloid_2021,
    title       = {Amyloid $\beta$ and {Amyloid Precursor Protein} synergistically suppress large-conductance calcium-activated potassium channel in cortical neurons},
    author      = {Yamamoto, K. and Yamamoto, R. and Kato, N.},
    journal     = {Frontiers in Aging Neuroscience},
    year        = {2021},
    volume      = {13},
    pages       = {660319},
    doi         = {10.3389/fnagi.2021.660319},
}

@article{zhang_effects_2014,
    title       = {Effects of amyloid $\beta$-peptide fragment 31–35 on the {BK} channel-mediated $\mathrm{K}^{+}$ current and intracellular free $\mathrm{Ca}^{2+}$ concentration of hippocampal {CA}1 neurons},
    author      = {Zhang, Y. and Shi, Z. and Wang, Z. and Li, J. and Chen, J. and Zhang, C.},
    journal     = {Neuroscience Letters},
    year        = {2014},
    volume      = {568},
    pages       = {72--76},
    doi         = {10.1016/j.neulet.2014.03.028}
}

@article{fischl_freesurfer_2012,
	title      = {{FreeSurfer}},
	journal    = {NeuroImage},
	author     = {Fischl, B.},
	year       = {2012},
	volume     = {62},
	number     = {2},
	doi        = {10.1016/j.neuroimage.2012.01.021},
	pages      = {774--781},
}

@article{jhanandas_cellular_2001,
    title       = {Cellular mechanisms for amyloid $\beta$-protein activation of rat cholinergic basal forebrain neurons},
    author      = {Jhanandas, J. H. and Cho, C. and Jassar, B. and Harris, K. and Tavish, D. M. and Easaw, J.},
    journal     = {Journal of Neurophysiology},
    year        = {2001},
    volume      = {86},
    number      = {3},
    pages       = {1067--1510},
    doi         = {10.1152/jn.2001.86.3.1312}
}

@article{good_beta-amyloid_1996,
    title       = {Beta-amyloid peptide blocks the fast-inactivating $\mathrm{K}^+$ current in rat hippocampal neurons},
    author      = {Good, A. T. and Smith, D. O. and Murphy, R. M.},
    journal     = {Biophysical Journal},
    year        = {1996},
    volume      = {70},
    pages       = {296--304},
    doi         = {10.1016/S0006-3495(96)79570-X},
}

@article{cangiani_hp-version_2014,
  title     = {hp-version discontinuous {G}alerkin methods on polygonal and polyhedral meshes},
  author    = {Cangiani, Andrea and Georgoulis, Emmanuil H and Houston, Paul},
  journal   = {Mathematical Models and Methods in Applied Sciences},
  volume    = {24},
  number    = {10},
  pages     = {2009--2041},
  year      = {2014},
}

@book{cangiani_hp-version_2017,
    title       = {hp-version discontinuous {G}alerkin methods on polygonal and polyhedral meshes},
    author      = {Cangiani, Andrea and Dong, Zhaonan and Georgoulis, Emmanuil H and Houston, Paul},
    year        = {2017},
    publisher   = {Springer},
}

@article{schwartz2016analytic,
    title   = {Analytic modeling of neural tissue: I. A spherical bidomain},
    author  = {Schwartz, Benjamin L and Chauhan, Munish and Sadleir, Rosalind J},
    journal = {The Journal of Mathematical Neuroscience},
    volume  = {6},
    number  = {1},
    pages   = {9},
    year    = {2016},
    doi     = {https://doi.org/10.1186/s13408-016-0041-1}
}

@article{antonietti_hp_2013,
    title       = {hp-version composite discontinuous {G}alerkin methods for elliptic problems on complicated domains},
    author      = {Antonietti, Paola F and Giani, Stefano and Houston, Paul},
    journal     = {SIAM Journal on Scientific Computing},
    volume      = {35},
    number      = {3},
    pages       = {A1417--A1439},
    year        = {2013},
}

@article{arnold_unified_2002,
    title       = {Unified analysis of discontinuous {G}alerkin methods for elliptic problems},
    author      = {Arnold, Douglas N and Brezzi, Franco and Cockburn, Bernardo and Marini, L Donatella},
    journal     = {SIAM journal on numerical analysis},
    volume      = {39},
    number      = {5},
    pages       = {1749--1779},
    year        = {2002},
}

@article{leimer_saglio_p-adaptive_2025,
    title       = {A p-adaptive polytopal discontinuous {G}alerkin method for high-order approximation of brain electrophysiology},
    author      = {Leimer Saglio, Caterina Beatrice and Pagani, Stefano and Antonietti, Paola F.},
    journal     = {Computer Methods in Applied Mechanics and Engineering},
    year        = {2025},
    volume      = {446},
    pages       = {118249},
    doi         = {10.1016/j.cma.2025.118249},
}

@misc{leimer_saglio_high-order_2024,
    title       = {A high-order discontinuous {G}alerkin method for the numerical modeling of epileptic seizures},
    author      = {Leimer Saglio, Caterina Beatrice and Pagani, Stefano and Corti, M. and Antonietti, Paola F.},
    eprint      = {ArXiv},
    doi         = {10.48550/arXiv.2401.14310},
    year        = {2024}
}

@article{antonietti_polytopal_2024,
    author      = {Antonietti, P. F. and Corti, M. and Martinelli, G.},
    title       = {Polytopal mesh agglomeration via geometrical deep learning for three-dimensional heterogeneous domains},
    journal     = {Mathematics and Computers in Simulation},
    year        = {2026},
    volume      = {241},
    number      = {Part B},
    pages       = {335--353},
    doi         = {10.1016/j.matcom.2025.10.019}
}

@article{yang_alzheimer_2022,
    title       = {Alzheimer’s disease and epilepsy: {A}n increasingly recognized comorbidity},
    author      = {Yang, F.  and Chen, L.  and Yu, Y.  and Xu, T.  and Chen, L.  and Yang, W.  and Wu, Q.  and Han, Y.},
    journal     = {Frontiers in Aging Neuroscience},
    year        = {2022},
    volume      = {14},
    doi         = {10.3389/fnagi.2022.940515},
}

@article{romoli_amyloid_2021,
    title       = {Amyloid-$\beta$: a potential link between epilepsy and cognitive decline},
    author      = {Romoli, M. and Sen, A. and Parnetti, L. and Calabresi, P. and Costa, C.},
    journal     = {Nature Reviews Neurology},
    year        = {2021},
    volume      = {17},
    pages       = {469--480},
    doi         = {10.1038/s41582-021-00505-9},
}

@article{antonietti_lymph_2025, 
    title       = {lymph: discontinuous po{LY}topal methods for {Multi-PHysics} differential problems}, 
    author      = {Antonietti, P. F. and Bonetti, S. and Botti, M. and Corti, M. and Fumagalli, I. and Mazzieri, I.}, 
    journal     = {ACM Transactions on Mathematical Software}, 
    year        = {2025}, 
    volume      = {51},
    number      = {1},
    pages       = {3:1–3:22},
    doi         = {10.1145/3716310}, 
}

@article{colli_franzone_spreading_1993, 
    title       = {Spreading of excitation in {3--D} models of the anisotropic cardiac tissue. {I. V}alidation of the eikonal model}, 
    author      = {Colli Franzone, P. and Guerri, L.}, 
    journal     = {Mathematical Biosciences}, 
    year        = {1993}, 
    volume      = {113},
    number      = {2},
    pages       = {145-209},
    doi         = {10.1016/0025-5564(93)90001-Q}, 
}

@article{zhang_clinical_2022, 
    title       = {The clinical correlation between {A}lzheimer's disease and epilepsy}, 
    author      = {Zhang, D. and Chen, S. and Xu, S. and Wu, J. and Zhuang, Y. and  Cao, W. and Chen, X. and Li, X.}, 
    journal     = {Frontiers in Neurology}, 
    year        = {2022}, 
    volume      = {13},
    pages       = {922535},
    doi         = {10.3389/fneur.2022.922535}, 
}

@article{bloom_amyloid_2014,
    title     = {Amyloid-$\beta$ and tau: the trigger and bullet in {A}lzheimer disease pathogenesis},
    author    = {Bloom, George S},
    journal   = {JAMA neurology},
    volume    = {71},
    number    = {4},
    pages     = {505--508},
    year      = {2014},
    doi       = {10.1001/jamaneurol.2013.5847},
}

@article{beghi_epidemiology_2020,
    title   = {The epidemiology of epilepsy},
    author  = {Beghi, Ettore},
    journal = {Neuroepidemiology},
    number  = {2},
    volume  = {54},
    year    = {2020},
    pages   = {185--191},
    doi     = {10.1159/000503831},
}

@article{hodgkin_components_1952,
    title       = {The components of membrane conductance in the giant axon of {Loligo}},
    author      = {Hodgkin, A. L. and Huxley, A. F.},
    journal     = {The Journal of Physiology},
    year        = {1952},
    volume      = {116},
    number      = {4},
    pages       = {473},
    doi         = {10.1113/jphysiol.1952.sp004718},
}

@article{hodgkin_quantitative_1952,
    title       = {A quantitative description of membrane current and its application to conduction and excitation in nerve},
    author      = {Hodgkin, A. L. and Huxley, A. F.},
    journal     = {The Journal of physiology},
    year        = {1952},
    volume      = {117},
    number      = {4},
    pages       = {500},
    doi         = {10.1113/jphysiol.1952.sp004764},
}

@article{quarteroni_integrated_2017,
    title       = {Integrated {Heart—Coupling} multiscale and multiphysics models for the simulation of the cardiac function},
    author      = {Quarteroni, A.M. and Lassila, T. and Rossi, S. and Ruiz-Baier, R.},
    journal     = {Computer Methods in Applied Mechanics and Engineering},
    year        = {2017},
    volume      = {314},
    pages       = {345-407},
    doi         = {10.1016/j.cma.2016.05.031},
}

@article{potse_comparision_2006,
    title       = {A comparison of monodomain and bidomain reaction-diffusion models for action potential propagation in the human heart},
    author      = {Potse, M. and Dube, B. and Richer, J. and Vinet, A. and Gulrajani, R.M.},
    journal     = {IEEE Transactions on Biomedical Engineering},
    year        = {2006},
    volume      = {53},
    issue       = {12},
    pages       = {2425-2435},
    doi         = {10.1109/TBME.2006.880875},
}

@article {corti_discontinuous_2023,
    title       = {Discontinuous {G}alerkin methods for {F}isher-{K}olmogorov equation with application to {$\alpha$}-synuclein spreading in {P}arkinson's disease},                                  
    author      = {Corti, M. and Bonizzoni, F. and Dede', L. and Quarteroni, A. M. and Antonietti, P. F.},
    journal     = {Computer Methods in Applied Mechanics and Engineering},
    volume      = {417},
    year        = {2023},
    pages       = {116450},
    doi         = {10.1016/j.cma.2023.116450},
}

@article{antonietti_discontinuous_2024, 
    author      = {Antonietti, P. F. and Bonizzoni, F. and Corti, M. and Dall'Olio, A.},
    title       = {{Discontinuous Galerkin} approximations of the heterodimer model for protein-protein interaction},
    journal     = {Computer Methods in Applied Mechanics and Engineering},
    year        = {2024},
    volume      = {431},
    pages       = {117282},
    doi         = {10.1016/j.cma.2024.117282},
}

@misc{lamontagne_oasis_2019,
    author       = {LaMontagne, P. J. and Benzinger, T. LS. and Morris, J. C. and Keefe, S. and Hornbeck, R. and Xiong, C. and Grant, E. and Hassenstab, J. and Moulder, K. and Vlassenko, A. G. and Raichle, M. E. and Cruchaga, C. and Marcus, D.},
    title      = {{OASIS-3: L}ongitudinal Neuroimaging, Clinical, and Cognitive Dataset for Normal Aging and {A}lzheimer Disease},
    year        = {2019},
    eprint      = {MedRxiv},
    doi         = {10.1101/2019.12.13.19014902},
}

@article{antonietti_agglomeration_2022,
    title       = {Agglomeration of polygonal grids using graph neural networks with applications to multigrid solvers},
    author      = {P. F. Antonietti and N. Farenga and E. Manuzzi and G. Martinelli and L. Saverio},
    journal     = {Computers \& Mathematics with Applications},
    year        = {2024},
    volume      = {154},
    pages       = {45--57},
    doi         = {10.1016/j.camwa.2023.11.015},
}
\end{document}